\documentclass[11pt]{amsart}

\usepackage{amscd}
\usepackage{amsmath}
\usepackage{amsthm}
\usepackage{amsfonts}
\usepackage{eucal}
\usepackage[all]{xy}

\newtheorem{theorem}{Theorem}[section]
\newtheorem{lemma}[theorem]{Lemma}
\newtheorem{proposition}[theorem]{Proposition}
\newtheorem{corollary}[theorem]{Corollary}

\theoremstyle{definition}

\newtheorem{definition}[theorem]{Definition}
\newtheorem{example}[theorem]{Example}
\newtheorem{subexamples}[theorem]{Subexamples}

\newtheorem{properties}[theorem]{Properties}
\newtheorem{facts}[theorem]{Facts}
\theoremstyle{remark}

\newtheorem{remark}[theorem]{Remark}

\newcommand{\field}[1]{\ensuremath{\mathbb{#1}}}

\newcommand{\C}{\field{C}}
\newcommand{\F}{\field{F}}
\newcommand{\N}{\field{N}}
\newcommand{\Pp}{\field{P}}

\newcommand{\T}{\field{T}}
\newcommand{\Z}{\field{Z}}

\newcommand{\calm}{\mathcal{M}}

\newcommand{\calo}{\mathcal{O}}

\newcommand{\bone}{\mathbf{1}}

\newcommand{\mba}{\mathbf{a}}

\newcommand{\lp}{\left( }
\newcommand{\rp}{\right) }

\newcommand{\gr}[2]{G_{#1}(\Pp^{#2})}
\newcommand{\fl}[2]{F(#1;\Pp^{#2})}
\newcommand{\py}{p_\infty}
\newcommand{\cq}{/\!\!/ }
\newcommand{\scho}[2]{\omega^{#2}_{#1}}
\newcommand{\sch}[2]{\Omega^{#2}_{#1}}

\newcommand{\atm}{\mbox{${\mathfrak{Atm}}$}}
\newcommand{\atmp}{\mbox{${\mathfrak{Atm}_p}$}}
\newcommand{\gar}[1]{\ensuremath{{\mathfrak{A}}_R\hspace{-.04cm}\langle #1 \rangle }}
\newcommand{\garf}[1]{\ensuremath{{\mathfrak{A}}_{R}^{fin}\hspace{-.06cm}\langle #1 \rangle }}
\newcommand{\cvp}[2]{\ensuremath{{\mathcal C}_{#1}\bigl(#2\bigr)}}
\newcommand{\cvpd}[3]{\ensuremath{{\mathcal C}_{#1,#2}\left( #3\right)}}
\newcommand{\apx}[1]{\ensuremath{{\Pi}_{p}(#1)}}

\newcommand{\ax}[2]{\ensuremath{\Pi_{#1}\bigl(#2\bigr)}}
\newcommand{\alg}[2]{\ensuremath{\mathcal{A}_{#1}\bigl(#2\bigr)}}
\newcommand{\algp}[2]{\ensuremath{\mathcal{A}^{\geq}_{#1}\bigl(#2\bigr)}}
\newcommand{\mmp}[2]{\ensuremath{M_{#1}\bigl(#2\bigr)}}
\newcommand{\hp}[3]{\ensuremath{H_{#1}(#2,#3)}}

\newcommand{\Zplus}{\ensuremath{{\field{Z}}_{+}}}

\newcommand{\ahcs}[2]{\ensuremath{P_{#1}^{#2}(t)}}

\newcommand{\aecs}[2]{\ensuremath{E_{#1}\bigl(#2\bigr)}}
\newcommand{\ecs}[2]{\ensuremath{E_{#1}\bigl(#2\bigr)}}
\newcommand{\pb}[1]{\ensuremath{{\Pp}\bigl(#1\bigr)}}
\newcommand{\vb}[2]{\ensuremath{{\mathcal O}_{#1}\bigl(#2\bigr)}}
\newcommand{\eone}{\ensuremath{E\oplus{\bone}}}
\newcommand{\eot}{\ensuremath{E_1\oplus E_2}}
\newcommand{\efib}{\ensuremath{{\pb{E_1}\!\!\times_{W}\!\pb{E_2}}}}

\newcommand{\eaeq}{\sim_{\text{alg}^+}}
\newcommand{\aeq}{\sim_{\text{alg}}}
\newcommand{\equdef}{:=}
\DeclareMathOperator{\supp}{supp}

\DeclareMathOperator{\Pic}{Pic}
\DeclareMathOperator{\rank}{rk}
\DeclareMathOperator{\effd}{Div_+}

\title[Chow quotients and Euler-Chow series]{Chow quotients and
projective bundle formulas for Euler-Chow series} 
\author{E. Javier Elizondo}
\address{Instituto de Matem\'aticas, UNAM, Mexico}
\thanks{The first author was supported in part by grants UNAM-DGAPA
IN101296 and CONACYT 3936-E} 
\author{P. Lima-Filho}
\address{Department of Mathematics, Texas A{\&}M University, USA}
\thanks{The second author was partially supported by NSF grant \#
DMS-9401533} 

\begin{document}
{\flushleft To appear in  {\it Journal of Algebraic Geometry}}
\begin{abstract}
Given a projective algebraic variety $X$, let $\Pi_p(X)$ denote the
monoid of {\em effective algebraic equivalence classes of effective
algebraic cycles} on $X$. The $p$-th Euler-Chow series of $X$
is an element in the formal monoid-ring  
$\Z[ \! [ \Pi_p(X) ]\! ]$
defined in terms of Euler characteristics of 
the Chow varieties $\cvpd{p}{\alpha}{X}$ of $X$, with 
$\alpha \in\Pi_p(X)$.
We provide a systematic treatment of such series,
and give  projective bundle formulas which generalize previous  
results by \cite{law&yau-hosy} and \cite{eli-tor}.
The techniques used involve the {\em Chow quotients} 
introduced in \cite{kap&stu&zel-quot}, and this allows the computation of various
examples including some Grassmannians and flag varieties.
There are relations between these examples and
representation theory, and further 
results point to interesting connections
between  Euler-Chow series for certain varieties
and the topology of the moduli spaces $\overline{\calm}_{0,n+1}$. 
\end{abstract}

\maketitle
\tableofcontents

\section{Introduction}
\label{sec:intro}

The use of topological invariants on moduli spaces has played a
vital role in various branches of mathematics and mathematical
physics in the last two decades. A light sampling under this
vast umbrella includes works in gauge theory, the theory of
instantons, various moduli spaces of vector bundles, moduli
spaces of curves and their compactifications, Chow varieties and
Hilbert schemes.

In this work we study a class of invariants for projective
varieties arising from the Euler characteristics of their Chow
varieties. These invariants first appeared in the work of 
H. B. Lawson and Steve S. T. Yau \cite{law&yau-hosy}, whose techniques
play an important role in this paper,
and they present, in various instances,  a quite nice and elegant
behavior which can often be codified in simple generating
functions. 

As a motivation, we start with some particular cases, which
are well studied in the literature. Let $X$ be a connected
projective variety and let $SP(X)$ denote the disjoint union
$\coprod_{d\geq 0} SP_d(X)$ of all symmetric products of $X$,
with the disjoint union topology, where $SP_0(X)$ is a single point.
One can define a function
$E_0(X) \,: \, \Z_+ = \pi_0(SP(X)) \to \Z$ which sends $d$ to
the Euler characteristic $\chi(SP_d(X))$ of the $d$-fold
symmetric product of $X$. This is what we call the {\em $0$-th
Euler-Chow function of $X$.} The same information
can be codified as a formal power series 
$E_0(X) = \sum_{d\geq 0} \chi(SP_d(X))t^d$, 
and a  result of Macdonald~\cite{mcd-sym} shows that $E_0(X)$
is given by the rational function  
$E_0 ( X ) =  \lp 1/(1-t) \rp^{\chi(X)}$. 

Another familiar instance arises in the case of divisors.
Given an $n$-dimensional projective variety $X$,
let $\effd(X)$ denote the space of effective divisors on $X$ and
assume that $\Pic_0(X) = \{0 \}$. Consider the function
$E : \Pic(X) \to \Z$ which sends $L\in \Pic(X)$ to $\dim
H^0(X,\calo(L))$. Observe that
\begin{enumerate}
\item Given $L\in \Pic(X)$, then $E(L) \neq 0$ if and only if 
$L = \calo(D)$ for some effective divisor $D$;
\item Under the given hypothesis, algebraic and linear
equivalence coincide, and two effective divisors $D$ and $D'$
are algebraically equivalent if and only if they are in the same
linear system.
\end{enumerate}
The last observation implies that $\effd(X)$ can be written as
$\effd(X) = \coprod_{\alpha \in \algp{n-1}{X}} 
\effd(X)_\alpha$, where $\algp{n-1}{X}$ 
is the monoid of algebraic equivalence classes of effective
divisors (cf. Fulton \cite[\S 12]{ful-inter}), and $\effd(X)_\alpha$ is the
linear system associated to $\alpha \in \algp{n-1}{X}$.
The first observation shows that the only relevant data 
to $E$ is given by $\algp{n-1}{X} \subset \Pic(X)$. Therefore, we
might as well restrict $E$ and define the {\em $(n-1)$-st
Euler-Chow function} of $X$ as the function
$\ecs{n-1}{X} : \algp{n-1}{X} \to \Z_+$ which sends
$\alpha \in \algp{n-1}{X}$ to the Euler characteristic
$\chi(\effd(X)_\alpha) = \dim H^0(X, \calo(L_\alpha))$, where
$L_\alpha$ is the line bundle associated to $\alpha$.

\begin{example}
\label{exmp:hfunction}
An even more restrictive case arises when $\Pic(X) \cong \Z$, and 
$\algp{n-1}{X} \cong \Z_+$ is generated by the class of a very ample
line bundle $L$. Then the $(n-1)$-st Euler-Chow function
$\ecs{n-1}{X} = \sum_{d\geq 0} \dim{H^0(X;\calo(L^{\otimes n}))}\  t^n$ 
is just the {\em Hilbert function} associated to the
projective embedding of $X$ induced by $L$. This is once again a
rational function.
\end{example}

In general,  the situation is not so simple, and we need to
introduce additional notions in order to approach cycles of arbitrary
dimension.  We start with the {\em Chow monoid} $\cvp{p}{X}$ of
effective $p$-cycles on $X$, which can be written as a disjoint
union $\ \coprod_{\alpha \in \ax{p}{X}} \cvpd{p}{\alpha}{X}\ $ of
connected projective (Chow) varieties $\cvpd{p}{\alpha}{X}$; cf.
Section \ref{sec:chow}. Here, $\ax{p}{X}= \pi_0(\cvp{p}{X})$
denotes the monoid of {\em {\bf\it{effective}} algebraic equivalence
classes of effective $p$-cycles}. This monoid should be
contrasted with $\algp{p}{X}$, the monoid of algebraic
equivalence classes of effective $p$-cycles; cf. 
\cite[\S 12]{ful-inter}.
In fact, there is a finite surjective monoid morphism
$\ax{p}{X} \to \algp{p}{X}$, and the Grothendieck group
associated to both monoids is $\alg{p}{X}$, the group of
algebraic equivalence classes of $p$-cycles on $X$; cf. Friedlander~\cite{fri-cycles}.
The properties and relations among these monoids is discussed in
Section~\ref{sec:chow}.

The {\em $p$-th Euler-Chow function} of the projective
variety $X$ is then defined as the function
\begin{align}
\label{eqn:ecs0}
\ecs{p}{X} : \ax{p}{X} & \longrightarrow \Z \\
\alpha  & \longmapsto \chi\lp \cvpd{p}{\alpha}{X} \rp \notag
\end{align}
which can be completely encoded as a formal power series
$\ecs{p}{X} = \sum_{\alpha \in \ax{p}{X}} \chi\lp
\cvpd{p}{\alpha}{X} \rp \ t^\alpha$,
on variables $t^\alpha$ associated to the elements
of $\ax{p}{X}$, and satisfying the relations
$t^{\alpha}t^{\beta}=t^{\alpha + \beta}$. For this reason, we also call
$\ecs{p}{X}$ the {\em $p$-th Euler-Chow series} of $X$.

Previous works in which such functions have been explicitly computed 
include the cases of products of two projective spaces in Lawson and Yau
\cite{law&yau-hosy}, simplicial projective toric varieties in
\cite{eli-tor}, and some cases of principal bundles whose structure
group is an abelian variety in Elizondo and Hain \cite{eli&hain-abe}.

The two main results of this paper, Theorems \ref{thm:main1} and
\ref{thm:cq} consider instances where 
one can reduce the calculation of Euler-Chow functions to simpler
situations. 
In the first case we consider the projectivization $\Pp(E)$ of an
algebraic vector bundle $E$ over a projective variety $W$, which 
splits as a direct sum of bundles $E_1\oplus E_2$. We define a ``trace
map'' $t_{p-1} : \cvp{p-1}{\Pp(E_1)\times_W  \Pp(E_2)} \to
\cvp{p}{\Pp(E)}$ which, when combined with the inclusions
$i_1 : \Pp(E_1) \to \Pp(E)$ and $i_2 : \Pp(E_2) \to \Pp(E)$ produces a
monoid morphism
$\Psi_p : \ax{p-1}{\Pp(E_1)\times_W \Pp(E_2)} \times  \ax{p}{\Pp(E_1)}
\times \ax{p}{\Pp(E_2)}\to \ax{p}{\Pp(E)}$. 
This morphism is the key ingredient for the 
following ``split bundle'' formula.
 
\noindent{\bf Theorem \ref{thm:main1}.}\ {\it 
Let $E_1$ and $E_2$ be algebraic vector bundles over a connected projective
variety $W$, of ranks $e_1$ and $e_2$, respectively, and let $0 \leq
p \leq e_1+e_2-1$. Then the $p$-th Euler-Chow function of $\ \pb{\eot}$ 
is given by
\begin{equation*}
\ecs{p}{\pb{\eot}} \, = \, 
{\Psi_{p}}_{\sharp} \lp\ \ecs{p-1}{\efib} \odot
\ecs{p}{\pb{E_1}} \odot \ecs{p}{\pb{E_2}}\  \rp.
\end{equation*}
}

\noindent In the formula, ${\Psi_p}_\sharp$ denotes the push-forward
induced by $\Psi_p$, introduced in Definition \ref{def:fund}, and
$\odot$ denotes the ``exterior product'' of Euler-Chow functions, cf. Definition
\ref{defn:exterior}. 

This formula applies to several situations, where the main technical
difficulty is reduced to  the computation of $\Psi_p$. 
For example, in the case where $E_1=L$ is
a line bundle generated by its global sections and $E_2$ is the trivial
line bundle, then $\Psi_p$ is completely determined in terms of the
first Chern class of $L$, as one expects; cf. Corollary \ref{cor:psi}.
In order to illustrate such a case, let $L$ be $ \calo_{\Pp^n}(d),\  d\geq 0$.
Then we have $\ax{p}{\Pp(\calo_{\Pp^n}(d)\oplus \bone )} \cong \Z_+\oplus
\Z_+$. If one writes $t^{(r,s)} = x^ry^s$,  then the $p$-th Euler-Chow
function of $\Pp\lp \calo_{\Pp^n}(d) \oplus \bone \rp $ 
is determined by the following generating function
\begin{equation*} 
\ecs{p}{\pb{\vb{\Pp^n}{d}\oplus \bone}} \, = \, \lp\frac{1}{1-x}\rp^{\lp\substack{n+1\\p}\rp}
 \lp\frac{1}{1-y}\rp^{\lp\substack{n+1\\p+1}\rp}\
 \lp\frac{1}{1-x^dy}\rp^{\lp\substack{n+1\\p+1}\rp}
\end{equation*}
cf. (\ref{2.21}).

\begin{remark}
The case $d=0$, i.e., $\Pp^n \times \Pp^1$ was computed in
\cite{law&yau-hosy}. The Euler-Chow series of the
blow-up $\tilde{\Pp}^{n+1}\cong \Pp(\calo_{\Pp^n}(1)\oplus \bone)$ 
of $\Pp^{n+1}$ at a point, and of the Hirzebruch surfaces $F_d \cong
\Pp(\calo_{\Pp^1}(d)\oplus \bone )$ were computed in Elizondo~\cite{eli-tor},
using general techniques for toric varieties.
\end{remark}

The next class of results applies to the case where $X$ is a smooth
projective variety on which $\C^*$ acts via automorphisms, 
and such that the
fixed point set $X^{\C^*}$  has only two connected components $X_1$ and
$X_2$. In this case, following Kapranov, Sturmfels and
Zelevinski~\cite{kap&stu&zel-quot}, we introduce the Chow 
quotient $X\cq \C^*$ of $X$ by $\C^*$, which comes equipped once again
with a ``trace map'' $t_{p-1} : \cvp{p-1}{X\cq {\C}^*}\to \cvp{p}{X}$, 
defined using techniques from  Friedlander and Lawson \cite{fri&law-cocyc}.

In a similar fashion to the previous theorem 
we introduce a monoid morphism
$\Psi_p : \ax{p-1}{X\cq \C^*} \times \ax{p}{X_1}\times \ax{p}{X_2} \to
\ax{p}{X}$
which produces the following result.

\noindent{\bf Theorem \ref{thm:cq}.}\ {\it
Let $X$ be an smooth projective  variety 
on which $\C^*$ acts algebraically.
If $X^{\C^*}$ is the union of two connected components $X_1$ and $ X_2$,
then for each $0\leq p\leq \dim{X}$ one has
\begin{equation*}
\ecs{p}{X} = {\Psi_p}_\sharp \lp \ecs{p-1}{X\cq \C^*}\odot
                            \ecs{p}{X_1} \odot \ecs{p}{X_2} \rp.
\end{equation*}
}
\begin{remark}
In the case were $W$ is smooth, then Theorem \ref{thm:main1} can be
shown to be a particular case of Theorem \ref{thm:cq}.
\end{remark}

Examples of such Chow quotients 
and resulting trace maps $t_p$'s and homomorphisms $\Psi_p$'s, are
computed  in Section \ref{subsec:exmp}.
There we consider the case where
$\C^*$ acts linearly on the last coordinate of $\C^{n+1}$, inducing an
action on all (partial) flag varieties $\fl{d_1,\ldots,d_r}{n}$, with
$0\leq d_1< \cdots < d_r \leq n$. Let $G_d(\Pp^n)$ be the space of all 
$d$-planes in $\Pp^n$, then we show that $G_d(\Pp^n)\cq \C^*\cong
\fl{d-1,d}{n-1}$, which implies the following ``almost'' recursive
formula 
$$
\ecs{p}{\gr{d}{n}} = \Psi_p\lp \ecs{p-1}{\fl{d-1,d}{n-1}}\odot
\ecs{p}{\gr{d}{n-1}}\odot \ecs{p}{\gr{p-1}{n-1}} \rp.
$$

As another example, we show that $\fl{0,1}{n}\cq \C^*$
is the blow-up of $\Pp^{n-1}\times \Pp^{n-1}$ along the diagonal.
In this context, this variety is naturally identified with $\Pp^{n-1}[2]$,  
the compactification of the configuration space of two distinct points
in $\Pp^{n-1}$, introduced by Fulton and MacPherson in \cite{ful&macP-config}.
The main formula in Theorem \ref{thm:cq} does not apply in this case,
since the fixed point set of the $\C^*$ action on $\fl{0,1}{n}$ has
three connected components. Nevertheless, this Chow quotient can still
be used  to obtain explicit computations in some cases; cf. Theorem
\ref{thm:f012div}, Corollary \ref{cor:f012all} and
Proposition \ref{prop:grass}.
In the tables below we exhibit the Euler-Chow series of
$\fl{0,1}{2}$ and $\gr{1}{3}$. 
We use the fact that for these spaces the
monoids $\ax{*}{-}$ are freely generated by the classes of Schubert
cycles, and indicate in the table the appropriate association
\{ \ variables \ \} $\leftrightarrow$ \{ Schubert classes \}.
$$
\begin{array}{|l|l|}\hline
\fl{0,1}{2}  & \text{Euler-Chow series} \\ \hline\hline
t \leftrightarrow \scho{0;0,1}{} & E_0 = \frac{1}{(1-t)^6} \\ \hline
r \leftrightarrow \scho{0;0,2}{},\ s \leftrightarrow \scho{1;1,2}{}  & 
E_1 = \frac{1}{(1-r)^{3}(1-s)^3(1-rs)^3} \\ \hline
x \leftrightarrow \scho{1;1,2}{},\ y \leftrightarrow \scho{0;0,2}{}  & 
E_2 = \frac{1-xy}{(1-x)^3(1-y)^3}
\\ \hline
\end{array}
\quad
\begin{array}{|l|l|}\hline
 G_1\Pp^3  & \text{Euler-Chow series} \\ \hline\hline
t \leftrightarrow \scho{0,1}{} & E_0 = \frac{1}{(1-t)^6} \\ \hline
s \leftrightarrow \scho{0,2}{} & E_1 = \frac{1}{(1-s)^{12}} \\ \hline
x \leftrightarrow \scho{0,3}{},y \leftrightarrow \scho{1,2}{}  & 
E_2 = \frac{1}{(1-x)^4(1-y)^4(1-xy)^3}
\\ \hline
z \leftrightarrow \scho{1,3}{} & E_3 = \frac{1+z}{(1-z)^5} \\ \hline
\end{array}
$$

There are several open questions suggested by the results we have computed so
far. For example, what is the precise 
relation between the Euler-Chow series of
generalized flag varieties and the combinatorics of such a variety? Even
in apparently simple cases, such as $\gr{1}{n}$ (not fully computed yet), one
finds a direct relation between their Euler-Chow series and the
ones for the compactified moduli space $\overline{\calm}_{0,n+1}$ of stable
$(n+1)$-punctured curves of genus $0$. The latter space arises as 
the Chow quotient $\gr{1}{n}\cq (\C^*)^{n}$; cf. \cite{kap-cq}. 
Other questions, such as general blow-up formulas and rationality 
of Euler-Chow series (not expected in general) are also quite 
challenging.

This paper is organized as follows. In Section \ref{sec:prelim}
we provide the necessary algebraic terminology and background. This
material is complemented with Appendix \ref{app:A}, where we describe
a more general approach to Euler-Chow series
using the notion of monoid-graded algebras and
their invariants. In Section \ref{sec:chow} we introduce the Chow
varieties and Chow monoids of projective varieties, along with various
properties. The monoids $\ax{p}{X}$ and $\algp{p}{X}$ are introduced and
compared  in this section, and relations between them and their common
Grothendieck group are presented. In this section we also introduce
the Euler-Chow functions.  In Section \ref{sec:proj} we study the case
of projective bundles and prove Theorem \ref{thm:main1}. We also study
the projective closure of line bundles and exhibit various explicit
examples. In Section \ref{sec:cq} we introduce the Chow quotient 
$X\cq (\C^*)^n$ of a projective variety $X$ under an algebraic action
of $(\C^*)^n$, following \cite{kap&stu&zel-quot}. We combine the
notions of Chow quotients and the trace maps of \cite{fri&law-cocyc}
to prove  Theorem \ref{thm:cq}. Various examples of Chow quotients, 
resulting trace maps, and Euler-Chow series are exhibited.
In the Appendix \ref{app:A} we present an algebraic framework which
places the Euler-Chow series $\ecs{p}{X}$
in a broader context, as an invariant of the Pontrjagin ring of the 
Chow monoid $\cvp{p}{X}$. Functoriality and ``change of monoid''
functors are studied in the more general context of monoids with
proper multiplication.
\medskip

\noindent{\em Acknowledgements:} The first author would like to thank both 
H. Blaine Lawson for suggesting a problem which gave origin to some of
the ideas in this article, and Alastair King for useful 
conversations about this topic.
The second author would like to thank
the {\em Instituto de Matem\'aticas, UNAM,} and the
{\em Sondersforschungsbereich 343/Osnabr\"uck Universit\"at} for their warm
hospitality and support during the elaboration of portions of this
work. He also thanks Roland Schw\"anzl for useful conversations on the
subject.

\section{Preliminaries}
\label{sec:prelim}

Let us start with an abelian monoid $M$, whose multiplication we denote by
$\ast_{M} : M \times M \longrightarrow M$.
When no confusion is likely to arise we use an additive notation 
$+ : M \times M \longrightarrow M$ with no subscripts attached.
We say that $M$ has {\bf finite multiplication} if $\ast_M$
has finite fibers. Typical examples are the freely generated monoids,
such as the non-negative integers $\Z_+$ under addition.

\begin{definition}
\label{defn:SM}
Given a monoid with finite multiplication $M$, and
a commutative ring $S$, denote by $S^M$ the set of all functions
from $M$ to $S$. If $f$ and $F^{\prime}$ are elements in $S^M$,
let $f + f^{\prime} \in S^M$ be defined
by pointwise addition, i.e. $(f + f^{\prime})(m) =
f(m) + f^{\prime}(m)$. \, Define the product 
$f \ast f^{\prime} \in S^M$ \, as the \lq\lq convolution\rq\rq\
$$
\lp f \ast f^{\prime}\rp (m) \, = \, 
\sum_{a\ast_M b = m} \, f(a) f^{\prime}(b).
$$
It is easy to see that $S^M$ then becomes a commutative ring
with unity, under these operations.
\end{definition}

\begin{remark}
\label{rem:power}
The ring $S^M$ can be identified with the completion $S[\![M]\!]$ of the
monoid algebra $S[M]$ at its augmentation ideal. Therefore, the
elements of $S^M$ can be written as a formal power series
$f= \sum_{m\in M} s_m\cdot t^m,$ on variables $t^m$ and coefficients
in $S$.  In this form the multiplication is given by the relation
$\ t^mt^{m'} = t^{m+m'}\ $ for elements $m,m'\in M$.
\end{remark}

\begin{definition}
\label{def:fund}
Given a monoid morphism $\Psi : M \longrightarrow N$,  
$f \in S^M$ and $g \in S^N$, \, define $\Psi^{\sharp}g \in S^M$ \, and \, 
$\Psi_{\sharp}f \in S^N$\, by
$$
\lp  \Psi^{\sharp}g\rp  (m) \, = \, g(\Psi (m))
$$
and
$$
\lp  \Psi_{\sharp}f\rp  (n) \,\, = 
 \sum_{m\in\Psi^{-1}(n)} \, f(m)
$$
if $\Psi$ has finite fibers.
\end{definition}
\begin{proposition}
\label{13}
Let $M$ and $N$ be monoids with finite multiplication,
and let $\Psi: M \longrightarrow N$ be
a monoid morphism. Then
\begin{enumerate}
\item The pull-back map
$\Psi^{\sharp} : S^N \longrightarrow S^M$ is an $S$-module 
homomorphism.

\item If $\Psi$ has finite fibers then 
the push-forward map
$\Psi_{\sharp} : S^M \longrightarrow S^N$ is a morphism of $S$-algebras.

\item Any ring homomorphism $\Psi : S \longrightarrow S^{\prime}$
induces a ring homomorphism 
$\Psi_{\ast} : S^M \longrightarrow {S^{\prime}}^{M}$.
\end{enumerate}
\end{proposition}

\begin{proof}\ 
{\bf 1.}\quad Given $f, g \in S^{M^\prime}$ and $m\in M, s\in S$, one has:
\begin{align*}
\Psi^{\sharp} (f+sg)(m) & = \lp   f+sg\rp   \lp  \Psi (m) \rp  
                    \; =\; f \lp  \Psi (m)\rp   + sg\lp   \Psi (m)\rp  \\ 
          & = \lp   \Psi^{\sharp} f\rp  (m) + s\lp  \Psi^{\sharp}g\rp  (m)
\end{align*}

{\bf 2.}\quad The proof that  $\Psi_{\sharp}$ is an $S$-module morphism 
follows the same pattern as the previos one. As to the multiplicative 
structure one has for $f,g \in S^M$ and  $n \in N$:
\begin{align*}
\Psi_{\sharp} \lp   f \ast g \rp  (n) 
        & = \sum_{m\in\Psi^{-1}(n)} 
                \lp   f \ast g \rp   (m)
         \; =\;  \sum_{m\in\Psi^{-1}(n)} \, \,
             \sum_{\substack{a,b\in{M}\\ a{\ast_{M}}b=m}} f(a) g(b)\\
        & = \sum_{\substack{\alpha ,\beta\in{N}\\\alpha\ast_{N}\beta=n}} 
           \, \, \sum_{\substack{a\in\Psi^{-1}(\alpha)\\b\in\Psi^{-1}(\beta)}}
                f(a) g(b)
       \; =\; \sum_{\substack{\alpha ,\beta\in{N}\\\alpha\ast_{N}\beta=n}}
            \lp   \Psi_{\sharp} f\rp   (a) \cdot \lp  \Psi_{\sharp} g\rp  (b)\\
        & = \lp  \Psi_{\sharp} f\rp   \ast \lp  \Psi_{\sharp} g\rp   (n) 
\end{align*}

{\bf 3.}\quad The proof is evident.
\end{proof}

The last operation we need to introduce is the following {\bf exterior
product}. 

\begin{definition}
\label{defn:exterior}
Given monoids $M$ and $N$, and a commutative ring $S$,
one can define a map $ \odot: S^M \otimes_S S^N  \; \rightarrow \;
S^{M\times N }$. This map sends $f\otimes g$ to the function
$f\odot g \; \in \; S^{M\times N }$\ which assigns
to $(m,n)$ the element $f(m)g(n) \in S.$
\end{definition}

\begin{proposition}\hfill
\label{prop:odot}
The operation $\odot$ is bilinear and associative. In other words,
the following diagram commutes:
$$
\xymatrix@C-=10pt{
&(S^M\otimes_S S^N)\otimes_S S^P \ar[rr]^-{\cong} 
\ar[dl]^{\odot \ \otimes \ 1} 
& & S^M\otimes_S (S^N\otimes_S S^P) 
\ar[dr]^{1 \ \otimes \ \odot}
& \\
S^{M\times N}\otimes_S S^P \ar[drr]_{\odot} & & 
& &  S^M\otimes S^{N\times P} \ar[dll]^{\odot} \\
&  & S^{M\times N\times P}
}
$$
\end{proposition}

\section{The Euler-Chow series of projective varieties}
\label{sec:chow}

Let $X$ be a projective algebraic variety over \C , and let $p$
be a non-negative integer such that $0 \leq p \leq \dim X$. 
The {\bf Chow monoid} \cvp{p}{X}\ of effective $p$-cycles on $X$ is the
free monoid generated by the irreducible $p$-dimensional 
subvarieties of $X$. It is well-known that $\cvp{p}{X}$ can be written
as a countable disjoint  union of projective algebraic varieties
$\cvpd{p}{\alpha}{X}$,  the so-called {\bf Chow varieties}. We
summarize, in the following statements, a few   basic properties  of
the  Chow monoids and varieties which are found in Hoyt~\cite{hoy-bun}, 
Friedlander~\cite{fri-cycles},and Friedlander and Mazur \cite{fri&maz-oper}. 
For a recent survey and
extensive bibliography on the subject, we refer the 
reader to Lawson~\cite{law-survey}.

\begin{properties} 
\label{proper:chow}
Let $X$  be a projective variety, and fix $0\leq p \leq \dim{X}$.
\begin{enumerate}
\item The disjoint union topology on $\cvp{p}{X}$ induced by the
classical topology on the Chow varieties, is independent of the
projective embedding of $X$; cf. \cite{hoy-bun}.
\item The restriction of the monoid addition to products of Chow
varieties is an algebraic map; \cite{fri-cycles}.
\item An algebraic map $f:X\rightarrow Y$ between projective
varieties (hence a proper map), 
induces a natural monoid morphism $f_* : \cvp{p}{X}
\rightarrow \cvp{p}{Y}$ which is an algebraic continuous map when
restricted to a Chow variety; cf. \cite{fri-cycles}. This is the {\em
proper push-forward map}.
\item A flat map $f:X\rightarrow Y$ of relative dimension $k$, 
induces a natural monoid morphism $f^* : \cvp{p}{Y}
\rightarrow \cvp{p+k}{X}$ which is an algebraic continuous map when
restricted to a Chow variety; cf. \cite{fri-cycles}. This is the {\em
flat pull-back map}.
\end{enumerate}  
\end{properties}
\begin{remark}
Since we work over $\C$, we can define an 
{\em algebraic continuous map} as a continuous map
$f: X\to Y$ between varieties which induces an algebraic map
$f^\nu : X^\nu \to Y^\nu$ between their weak normalizations; 
cf. \cite{fri&maz-oper}. One could alternatively 
define {\em Chow varieties} as the weak normalization of those we 
consider here, as in Koll\'ar~\cite{kol-rat}. This does not alter their topology,
but transforms Chow varieties into a functor in the category of
projective varieties and algebraic morphisms. Either approach can be
used in this paper, without altering the results.
\end{remark}

\begin{definition} \hfill \\
\noindent{\bf 1.}\quad We denote by \apx{X}, the monoid 
$\pi_0 (\cvp{p}{X})$ of path-components of \cvp{p}{X}. This is the monoid
of \, \lq\lq effective algebraic equivalence classes\rq\rq\ of effective
$p$-cycles on $X$. We use the notation $a\eaeq b$ to express
that two effective  cycles $a,b$ are effectively
algebraically equivalent.  

\noindent{\bf 2.}\quad The group of all algebraic $p$-cycles on $X$ modulo
algebraic equivalence is denoted $\alg{p}{X}$, and the submonoid of
$\alg{p}{X}$ generated by the classes of cycles with non-negative 
coefficients is denoted by $\algp{p}{X}$; cf. Fulton \cite[\S 12]{ful-inter}.  
We use the notation $a\aeq b$ to express that two cycles $a,b$ are 
algebraically equivalent.  

\noindent{\bf 3.}\quad Let $c:\cvp{p}{X}\ \longrightarrow \hp{2p}{X}{\Z}$
be the cycle map into singular homology; cf.~\cite[\S 19]{ful-inter}.
The image of $c$ is denoted by \mmp{p}{X}.
\end{definition}

The following result explains the relation between the monoids above.

\begin{proposition}\hfill 
\label{prop:gcp}        
\begin{enumerate} 
\item The Grothendieck group associated to the monoid $\ax{p}{X}$ is
$\alg{p}{X}$. In particular, there is a natural monoid morphism  
$\iota_p : \ax{p}{X} \rightarrow \alg{p}{X}$
which satisfies the universal property that any  monoid morphism 
$f: \ax{p}{X} \rightarrow G$, from $\ax{p}{X}$ into a group $G$,
factors through $\alg{p}{X}$.
\item The image of $\iota_p$ is $\; \algp{p}{X}$, and the image of
$\algp{p}{X}$ under the cycle map is $\mmp{p}{X}$.
\item The monoid surjection 
$\overline{\iota}_p : \ax{p}{X} \rightarrow \algp{p}{X}$ induced by
$\iota_p$ is an isomorphism if and only if $\ax{p}{X}$ has cancellation
law. 
\item Both $\overline{\iota}_p : \ax{p}{X} \to \algp{p}{X}$ and 
$c_p : \algp{p}{X} \to \mmp{p}{X}$ are finite monoid morphisms.
\end{enumerate}
\end{proposition}
\begin{proof}
The first assertion is proven in \cite{fri-cycles},
and follows from standard arguments, e.g. in  Samuel~\cite{sam-sem}. The second
assertion follows from the definitions and the universal property just
described. The third assertion follows from the elementary 
fact that the universal map from an abelian monoid into its group 
completion is injective if and only if the monoid has cancellation law. 
To prove the last assertion, consider a projective
embedding of $X$. The cycles supported in $X$ with
a fixed degree $d$ in the ambient projective space form a projective
variety, which is then a finite union of components of $\cvp{p}{X}$.
The assertion now follows easily from these observations.
\end{proof}

The following result is found in  Elizondo~\cite{eli-tor}.
\begin{proposition}
Given a complex projective algebraic variety $X$ and $0\leq p \leq \dim X$, 
the monoids \cvp{p}{X}, \ax{p}{X}, \algp{p}{X}\ and \mmp{p}{X}\ all have
finite multiplication. In particular, they all belong to \atmp; cf.
Appendix \ref{app:A}.
\end{proposition}

\begin{proof}
One has  surjective monoid morphisms:
\begin{equation}
\label{morphism}
\cvp{p}{X} \xrightarrow{\pi_p} \apx{X}
\xrightarrow{\iota_p} \algp{p}{X} \xrightarrow{c_p} \mmp{p}{X},
\end{equation}
so that when  a projective embedding $j: X \hookrightarrow \Pp^n$ is chosen
one obtains a commutative diagram
\begin{center}
\mbox{
\xymatrix{
\cvp{p}{X}  \ar[r] \ar[d]^{j_{\ast}}   &   \apx{X} \ar[r] \ar[d]^{j_{\ast}} & 
\algp{p}{X} \ar[r] \ar[d]^{j_{\ast}}   &   \mmp{p}{X} \ar[d]^{j_{\ast}}\\
\cvp{p}{\Pp^n} \ar[r]                  &   \apx{\Pp^n} \ar[r]^\cong \ar@{=}[d]  &  
\algp{p}{\Pp^n} \ar[r]^\cong \ar@{=}[d]& \mmp{p}{\Pp^n} \ar@{=}[d] \\
                                       & \Zplus & \Zplus & \Zplus
} }
\end{center}
where the leftmost vertical arrow is a closed inclusion. Recall that 
$\cvp{p}{\Pp^n} \, = \, \coprod_{d\in\Zplus} \cvpd{p}{d}{\Pp^n}$ 
where \cvpd{p}{d}{\Pp^n}\
is a projective connected algebraic (Chow) variety. Furthermore,
$\; j\lp  \cvp{p}{X}\rp  \cap \cvpd{p}{d}{\Pp^n}$ 
is a subvariety for all $d$.
It follows that \apx{X}, \algp{p}{X}\ and \mmp{p}{X}\ are all 
discrete, and that \cvp{p}{X}\ has 
finite multiplication, since it is free.
Proposition~\ref{pull-backs} now shows that \ax{p}{X}, \algp{p}{X}\  and
\mmp{p}{X}\  also have finite multiplication.
\end{proof}

\begin{definition}
\label{def:EC}
The (algebraic) {\bf $p$-th Euler-Chow function} of $X$ is the function
\begin{align}
\ecs{p}{X} : \ax{p}{X} & \longrightarrow \Z \\
\alpha &\longmapsto \chi(\cvpd{p}{\alpha}{X}), \notag
\end{align} 
which sends $\alpha \in \ax{p}{X}$ to
the topological Euler characteristic of $\cvpd{p}{\alpha}{X}$ (in the
classical topology).
\end{definition}

Following Remark \ref{rem:power} we associate a variable $t^\alpha$ to
each $\alpha \in \ax{p}{X}$ and express the $p$-th Euler-Chow function
as a formal power series
\begin{equation}
\label{eqn:power}
\ecs{p}{X}\  = \sum_{\alpha \in \ax{p}{X}} \chi\lp \cvpd{p}{\alpha}{X} \rp
\ t^\alpha.
\end{equation}

\begin{remark}\hfill
\label{rem:Pontrjagin}

\noindent{\bf 1.}\ The Pontrjagin ring 
\hp{\ast}{\cvp{p}{X}}{\Z}\ comes equipped with lots
of additional structure, which can be better expressed using the
terminology of Appendix \ref{app:A}.
It follows from Example~\ref{1.6} that it
becomes a finite \ax{p}{X}-graded ring. 
In other words, $\hp{\ast}{\cvp{p}{X}}{\Z} \in \garf{\ax{p}{X}}$.
Under this framework, given $X$ and $p$ as above, one could define 
the {\em (algebraic) Hilbert-Chow function}\ of $X$ as the function 
$\ahcs{p}{X} \in \Z[t]^{\ax{p}{X}}$ obtained as  
the {\em Hilbert \ax{p}{X}-series} of the Pontrjagin ring 
\hp{\ast}{\cvp{p}{X}}{\Z}\ in the sense of Definition ~\ref{defn:hilbert}.
In particular, the {\em algebraic Euler-Chow function}
$\aecs{p}{X} \in \Z^{\apx{X}}$ is the {\bf Euler \ax{p}{X}-series} of \ 
\hp{\ast}{\cvp{p}{X}}{\Z}; cf. Definition \ref{defn:hilbert}.

\noindent{\bf 2.}\ One could in a similar fashion define the 
$p$-th Euler-Chow function mapping either $\algp{p}{X}$ or 
$\mmp{p}{X}$ to $\Z$. These would simply be the functions 
$\overline{\iota}_\sharp(\ecs{p}{X})$ and 
$(c_p\circ\overline{\iota}_p)_\sharp (\ecs{p}{X})$; 
cf. Definition~\ref{def:fund}. 
\end{remark}
\begin{example}\hfill 

\noindent{\bf 1.}\ If $X$ is a connected variety, then $\cvp{0}{X} =
\coprod_{d\in \Z_+} SP_d(X)$, where $SP_d(X)$ is the $d$-fold symmetric
product of $X$. Therefore, the $0$-th Euler-Chow function is given by 
$$\ecs{0}{X} = \sum_{d\geq 0} \chi( SP_d(X) ) \ t^d = \lp \frac{1}{1-t}
\rp^{\chi(X)},$$
according to McDonald's formula \cite{mcd-sym}.

\noindent{\bf 2.}\  For $X = \Pp^n$, one has $\ax{p}{\Pp^n} \cong \Z_+$,
with the isomorphism given by the degree of the cycles. In this case, the
$p$-th Euler-Chow function was computed in \cite{law&yau-hosy}:
$$\ecs{p}{\Pp^n} = \sum_{d\geq 0} \chi(\cvpd{p}{d}{\Pp^n} ) \ t^d = 
\lp \frac{1}{1-t} \rp^{\binom{n+1}{p+1}}.$$

\noindent{\bf 3.}\ Suppose that $X$ is an $n$-dimensional variety such that
$Pic(X) \cong \Z$, generated by a very ample line bundle $L$. Then, we
have seen in Example \ref{exmp:hfunction} that 
$\ax{n-1}{X} \cong \Z_+$ and that
$\ecs{n-1}{X}$ is precisely the Hilbert function  
for the projective embedding  of $X$ induced by $L$.
\end{example}

\section{Projective bundle formulas}
\label{sec:proj}

In this section we exhibit a formula for 
the Euler-Chow function of certain projective
bundles over a variety $W$, and compute several examples.
The basic setup is the following. Consider two algebraic vector bundles
$E_1 \rightarrow W$ and $E_2 \rightarrow W$ over a complex projective 
variety $W$. The various maps involved in our discussion are displayed
in the commutative diagram below:
\begin{equation}
\label{eqn:notation}
\xymatrix{
\pb{E_1} \; \ar@{^{(}->}^-{i_1}[r] \ar[dr]_-{p_1} &
        \pb{E_1\oplus E_2} \ar[d]_-{q}&\
        \pb{E_2} \ar@{_{(}->}_-{i_2}[l] \ar[dl]^-{p_2}\\
& W
}
\end{equation}
where $p_1, p_2$ and $q$ are projections from the indicated
projective bundles, and
$i_1,i_2$ are the canonical inclusions. 

We will introduce a monoid monomorphim
$t_p : \cvp{p-1}{\efib} 
\longrightarrow \cvp{p}{\pb{\eot}}$ in Definition \ref{def:map1},
which is a closed inclusion. In this way we become equipped with
three morphisms of {\em monoids with finite multiplication}:
\begin{equation}
{i_1}_p : \ax{p}{\pb{E_1}} \longrightarrow \ax{p}{\pb{\eot}}
\end{equation}
induced by $ {i_1}$,
\begin{equation}
{i_2}_p : \ax{p}{E_2} \longrightarrow \ax{p}{\pb{\eot}}
\end{equation}
induced by ${i_2}$, and
\begin{equation}
\varphi_p : \ax{p-1}{\efib} \longrightarrow 
\ax{p}{\pb{\eot}}
\end{equation}
induced by $t_p$.

These three maps induce a morphism (with finite fibers)
\begin{equation}
\Psi_p : \ax{p-1}{\efib} 
\times \ax{p}{\pb{E_1}} \times \ax{p}{\pb{E_2}}
\longrightarrow \ax{p}{\pb{\eot}},
\end{equation}
by sending $(a,b,c)$ to $\Psi_p(a,b,c)= 
\varphi_p(a)+{i_1}_p(b) + {i_2}_p(c)$.

The main result in this section is the following.

\begin{theorem}
\label{thm:main1}
Let $E_1$ and $E_2$ be algebraic vector bundles over a connected projective
variety $W$, of ranks $e_1$ and $e_2$, respectively, and let $0 \leq
p \leq e_1+e_2-1$. Then the $p$-th Euler-Chow function of $\ \pb{\eot}$ 
is given by
\begin{equation}
\label{eqn:mformula}
\ecs{p}{\pb{\eot}} \, = \, 
{\Psi_{p}}_{\sharp} \lp\ \ecs{p-1}{\efib} \odot
\ecs{p}{\pb{E_1}} \odot \ecs{p}{\pb{E_2}}\  \rp.
\end{equation}
\end{theorem}

\begin{remark}
In case $\ \dim{X} < p$, the Chow monoid $\cvp{p}{X}$ consists of the zero
element only, which is the cycle with empty support. Therefore
$\ax{p}{X} = \{ 0 \}$ and $\ecs{p}{X} \equiv 1 \in \Z^{ \{ 0 \} }\equiv
\Z.$ 
\end{remark}

Before proving the theorem and providing examples, we must
define the map $t_p$ appearing in the formulas above. Let $L_1$
and $L_2$ denote the tautological line bundles $ \calo_{E_1}(-1)$ and $
\calo_{E_2}(-1)$ over $\pb{E_1}$ and $\pb{E_2}$, respectively, and let
$\pi_1$ and $\pi_2$ denote the respective projections from $\efib$
onto $\pb{E_1}$ and $\pb{E_2}$. The $\mathbb{P}^1$-bundle 
$\pi : \pb{\pi_1^*(L_1) \oplus \pi_2^* (L_2) } \rightarrow  \efib$
is precisely the blow-up of $\pb{\eot}$ along $\pb{E_1}\amalg \pb{E_2}$,
which we denote by $Q$, for short; see Lascou and Scott~\cite{las&sco-blow} for details. 
Let $b :  Q \rightarrow \pb{\eot}$ denote the blow-up map.

Since $\pi$ is a flat map of relative dimension $1$, and $b$
is a proper map, one has two algebraic continuous homomorphisms
(cf. \ref{proper:chow}), given by the flat pull-back
\begin{equation}
\pi^* : \cvp{p-1}{\efib} \rightarrow \cvp{p}{Q}
\end{equation}
and the proper push-forward
\begin{equation}
b_* : \cvp{p}{Q} \rightarrow \cvp{p}{\pb{\eot}}.
\end{equation}

\begin{definition}
\label{def:map1}
The map $t_p : \cvp{p-1}{\efib} \rightarrow \cvp{p}{\pb{\eot}}$ is
defined as the composition $ t_p = b_* \circ \pi^*$.
\end{definition}

The following lemma will be used as a reference in a few places throughout 
the paper.

\begin{lemma}
\label{lem:tech}
Let $X = \coprod_{i\in \N} X_i$ and $Y=\coprod_{j\in \N} Y_j $ be spaces
which are a countable disjoint union of connected projective varieties,
and let $f: X \rightarrow Y$ be a continuous map such that the
restriction  $f_{|_{X_i}}$ is an algebraic continuous map from $X_i$
into some $Y_j$. If $f$ is a bijection, then it is a homeomorphism in  
the classical topology.
\end{lemma}
\begin{proof}
Given $j\in \N$, one can write $f^{-1}(Y_j)$ as a disjoint union
$\coprod_{k \in \Lambda} X_{i_k}$, for a certain collection of connected
components of $X$. On the other hand, since $f$
is an algebraic continuous map, then $f(X_{i_k})$ is a Zariski closed
subset  of $Y_j$, and since $f$  is  one-to-one, then $Y_j$ 
is written as a union of a (at most countable)
family of disjoint closed algebraic subsets. This contradicts the
connectedness of $Y_j$, unless there is only one $j(i)$ such that
$f^{-1}(Y_{j})= X_{j(i)}$, in which case $f$ sends $X_{j(i)}$
homeomorphically onto $Y_j$.
\end{proof}

We now prove Theorem \ref{thm:main1}.
\begin{proof}[Proof of theorem \ref{thm:main1} ]
Consider the action of $\C^*$ on $\pb{\eot}$ given by
scalar multiplication on one of the factors of $\eot$, whose fixed point
set $\pb{\eot}^{\C^*}$ consists of $\pb{E_1}\amalg \pb{E_2}$. It is
easy to see that this induces an algebraic continuous 
action on $\cvp{p}{\pb{\eot}}$ via monoid automorphisms, 
and our next step is to identify its fixed point set.

Consider the map
\begin{equation}
\label{eqn:homeo}
\phi_p : \cvp{p-1}{\efib} 
\times \cvp{p}{\pb{E_1}} \times \cvp{p}{\pb{E_2}}
\longrightarrow \cvp{p}{\pb{\eot}},
\end{equation}
defined as $\phi_p(a,b,c) = t_p(a)+ i_{1*} (b) + i_{2*}(c)$. We claim
that $\phi_p$ is a homeomorphism onto the fixed point set
$\cvp{p}{\pb{\eot}}^{\C^*}$. 

If an element $\sigma = \sum_i n_i V_i \in \cvp{p}{\pb{\eot}}$ 
is fixed under the action, then  each  of its irreducible components
must clearly be invariant under the action, as well. On the other hand,
an invariant irreducible subvariety $V\subset \pb{\eot}$ 
may be of two types: 
\begin{description}
\item[i] those whose general
points are fixed under the action (and hence all points of $V$ are
fixed); 
\item[ii] and those whose general points have non-trivial orbits (of
dimension $1$).  
\end{description}
Recall that the Chow monoids $\cvp{p}{X}$  of any variety $X$  
are freely  generated by the irreducible $p$-dimensional subvarieties   
of $X$. 
Given cycles $\sigma_k \in \cvp{p}{\pb{E_k}}$, $k=1,2$, the support of
$i_{k*}(\sigma_k)$ is contained in $\pb{E_i}$, and hence
the images of $\cvp{p}{\pb{E_1}}$ and $\cvp{p}{\pb{E_2}}$ 
under $i_1$ and $i_2$
are freely generated by disjoint subsets of the 
generating set of $\cvp{p}{\pb{\eot}}$, consisting of varieties of the
first type described above.

On the other hand, given a $(p-1)$-dimensional subvariety $Z$ of
$\efib $, its inverse image $\pi^{-1}(Z)$ is a
$p$-dimensional subvariety of $Q$ whose points outside the exceptional
divisor of the blow-up map $b$,  have orbits of dimension $1$. Since
$b$ is a $\C^*$-equivariant birational map, the image $b(\pi^{-1}(Z))$
is an invariant irreducible subvariety of $\pb{\eot}$ of the second type.
It follows that the images of $i_{1*}, i_{2*}$ and $t_p$ are respectively
freely generated by disjoint subsets, and this proves that the map
$\phi_p$ is injective.

In order to show surjectivity, one just needs to show that every
invariant, irreducible subvariety $V \subset \pb{\eot}$ of the second type,
is of the form $b(\pi^{-1}(Z))$ for some $(p-1)$-dimensional subvariety
of $\efib$, since those of the first type are, by definition, in the
image of $i_{1*}$ and $i_{2*}$. Indeed, let $\tilde{V} \subset Q$ be the
proper transform of $V$ under $b$, and let $Z= \pi(\tilde{V}) \subset
\efib$. Since the general points of $V$ have orbits of dimension $1$, so
do the general points of $\tilde{V}$, and this shows that the general fiber
of $\pi_{|_{\tilde{V}}} : \tilde{V} \rightarrow Z$ has dimension $1$. In
particular, $\tilde{V} = \pi^{-1}(Z)$ and,
since  $\pi_{|_{\tilde{V}}} $  is a bundle projection and
$b$  sends $\tilde{V}$ birationally onto $V$, one concludes that
$t_p(Z) = b_*\circ\pi^*(Z)=V$.

The arguments above show that $\phi_p$ provides an algebraic continuous
bijection from the product $\cvp{p-1}{\efib} \times~\cvp{p}{\pb{E_1}}~\times~
\cvp{p}{\pb{E_2}}$ onto the fixed point set 
$\cvp{p}{\pb{E_1\oplus E_2}}^{\C^*}.\ $ 
Notice that for each $\alpha \in \ax{p}{\pb{\eot}}$ the fixed point set
$\cvpd{p}{\alpha}{\pb{\eot}}^{\C^*}$ is an algebraic subset of 
$\cvpd{p}{\alpha}{\pb{\eot}}$  (not necessarily connected), and hence
it can also be written as a countable disjoint union of projective
varieties.  One then applies Lemma \ref{lem:tech} to conclude that
$\phi_p$ is a homeomorphism onto its image. 

It is a general fact that the Euler characteristic of a variety with an
algebraic torus action equals that of its fixed point set; see 
Lawson and Yau \cite{law&yau-hosy} or Elizondo and Hain \cite{eli&hain-abe}.
Therefore, given $\alpha \in \ax{p}{\pb{E_1\oplus E_2}}$, it follows  that
\begin{equation}
\label{eqn:efix}
\chi\left( \cvpd{p}{\alpha}{\pb{E_1\oplus E_2}}  \right) =
\chi\left( \cvpd{p}{\alpha}{\pb{E_1\oplus E_2}}^{\C^*}  \right).
\end{equation}
On the other hand, if $\Psi_p$ is the morphism induced by $\phi_p$
between the monoids of  connected components 
then
\begin{multline}
\label{eqn:comp}
\cvpd{p}{\alpha}{\pb{E_1\oplus E_2}}^{\C^*} = \\
\coprod_{(a,b,c)\in \Psi_p^{-1}(\alpha)} 
\cvpd{p-1}{a}{\efib} \times \cvpd{p}{b}{\pb{E_1}}\times
\cvpd{p}{c}{\pb{E_2}}  ,
\end{multline}
since $\phi_p$ is a homeomorphism onto its image.
Therefore  (\ref{eqn:efix}) implies that
\begin{multline}
\chi \left( \cvpd{p}{\alpha}{\pb{E_1\oplus E_2}} \right) 
=\\
\sum_{(a,b,c)\in \Psi_p^{-1}(\alpha)} 
        \chi (\cvpd{p-1}{a}{\efib}  )\cdot
        \chi (\cvpd{p}{b}{\pb{E_1}} )\cdot
        \chi (\cvpd{p}{c}{\pb{E_2}} ) .
\end{multline}
It follows from Definitions \ref{def:map1} and \ref{def:fund} that
(\ref{eqn:mformula}) holds.
\end{proof}

\begin{remark}
The preceding proof also shows that the map $t_p$, 
introduced in Definition \ref{def:map1}, is a closed inclusion, as claimed. 
\end{remark}

\subsection{Projective closure of line bundles}

Here we  consider the case where $E_1 = E$ is a line bundle over a
projective variety $W$, and $E_2 = \bone_W$ is the trivial line bundle.
In this case, 
$$\Pp(E_1) = \Pp(E_2) = \Pp(E_1)\times_W \Pp(E_2) = W ,$$
and the inclusions 
$i_k : \Pp(E_k) \hookrightarrow \Pp(E_1\oplus E_2), \ k=1,2$ 
(cf. Diagram \ref{eqn:notation}), become sections  
$i : W \rightarrow \pb{\eone}$ and 
$j : W \rightarrow \pb{\eone}$ 
of the projective bundle $\pb{\eone}$ over $W$.

If $\xi = c_1\lp \vb{\eone}{1}\rp $ denotes the first Chern class of the
canonical bundle over \pb{\eone}, then the map
\begin{align}
\label{eqn:isomorphism}
T : \alg{p-1}{W} \oplus \alg{p}{W} & \rightarrow \alg{p}{\pb{\eone}} \\
(\alpha,\beta) & \longmapsto q^*\alpha + \xi \cap q^*\beta \notag
\end{align}
is an isomorphism. Note that, since $\xi \cap q^*\beta = i_* \beta$, the
isomorphism above restricts to an injection
\begin{equation}
\label{eqn:injection}
T^\geq : \algp{p-1}{W} \oplus \algp{p}{W} \rightarrow
\algp{p}{\pb{\eone}}.
\end{equation}

\begin{lemma}
\label{lem:sections}
Let $E$ be a line bundle over $W$ which is generated by its global sections. 
Then:
\begin{description}
\item[a] The injection $T^\geq$ is an isomorphism;
\item[b] If $\ax{*}{W}$ are monoids with cancellation law for every $*$,
then so are $\ax{*}{\pb{\eone}}$. Equivalently, if the natural
surjections $\ax{p}{W} \rightarrow \algp{p}{W}$ are isomorphisms for
all $p$, then so are the surjections
$\ax{p}{\pb{\eone}} \rightarrow \algp{p}{\pb{\eone}}$.
\end{description}
\end{lemma}
\begin{proof}
{\bf a.}
It follows from the proof of Theorem \ref{thm:main1} that every
effective $p$-cycle $a$ in $\pb{\eone}$ is effectively algebraically 
equivalent to a cycle of the form $ i_*b + j_* c + q^*d$, where
$b,c \in \cvp{p}{W}$ and $d \in \cvp{p-1}{W}$. Since $E$ is generated by
its sections, one can find a section  $s : W \rightarrow E$ 
of $E$ whose zero locus $Z\subset W$ intersects $c$ properly. Let
$\tilde{s} : X \rightarrow \pb{\eone}$ denote the composition 
$\iota \circ s$, where $\iota : E \hookrightarrow \pb{\eone}$ is the
open inclusion.  Then, the closure of the orbit of 
$\tilde{s}_* c$ under the $\C^*$
action on $\cvp{p}{\pb{\eone}}$ contains two fixed points: $j_* c$ and 
$i_* c + q^*(Z\cap c)$. It follows that 
$a \eaeq i_*(b+c) + q^*(d + Z\cap c)$, and this
shows that $T^\geq$ is surjective.

{\bf b.} We need to show that if $a\aeq a'$ then $a\eaeq a'$, where 
$a,a' \in \cvp{p}{\pb{\eone}}$. The arguments described in the first
half of the proof show that one can find cycles $b,b' \in \cvp{p}{W}$
and $c,c' \in \cvp{p-1}{W}$ such that $a\eaeq i_*b +q^*c$ and 
$a'\eaeq i_*b' +q^*c'$. Since $a\aeq a'$, one concludes that $b=q_*a
\aeq q_*a'=b'$, and hence the hypothesis implies that $b\eaeq b'$.
Also, $j_*W$ intersects properly both $i_*b+q^*c$ and $i_*b'+q^*c'$, 
and hence 
$j_*c = j_*W \cdot (i_*b+q^*c) \aeq j_*W \cdot ( i_*b'+q^*c') = j_*c'$.
Applying $q_*$ shows that $c \aeq c'$ and hence $c\eaeq c'$, and this
concludes the proof. The equivalence of both statements follows from
Proposition \ref{prop:gcp}.
\end{proof}

\begin{corollary}
\label{cor:psi}
Under the same hypothesis,  the homomorphism
$$\Psi_p : \ax{p-1}{W}\oplus \ax{p}{W}\oplus \ax{p}{W} \rightarrow
\ax{p}{\pb{\eone}} \cong \ax{p-1}{W} \oplus \ax{p}{W}$$
sends $(\alpha, \beta,\gamma)$ to $(\alpha+\beta,\ \xi\cap \gamma \ + \alpha)$.
\end{corollary}
\begin{proof}
This is clear.
\end{proof}

\begin{example}
Let $W = \Pp^n$ and $E = \vb{\Pp^n}{d} $, with $d \geq 0$. 
We will compute 
$\ecs{p}{ \pb{ \vb{\Pp^n}{d} \oplus \bone} }$ 
for $1 \leq p \leq n$. 
First observe that $W$ and $E$ satisfy the hypothesis of Lemma
\ref{lem:sections}, and that 
$$
\ax{p-1}{\pb{E}} \cong \ax{p}{\Pp^{n}}
 = \Zplus \cdot [\Pp^{p-1}] \cong \Zplus
\quad
\text{ and }
\quad
\ax{p}{\pb{E}} \cong \ax{p}{\Pp^n} \cong \Zplus \cdot [\Pp^p] 
\cong \Zplus ,
$$ 
where $[\Pp^{p}]$ denotes the  class of a $p$-dimensional linear subspace 
of $\Pp^n$. By definition, the map $\Psi_p$ fits into a 
commutative diagram
$$
\xymatrix{
 \ax{p-1}{\pb{\vb{\Pp^n}{d}}}  \times
\ax{p}{\pb{\vb{\Pp^n}{d}}} \times
\ax{p}{\Pp^n} \ar[r]^-{\Psi_p} \ar@{=}[d]
& \ax{p}{\pb{\vb{\Pp^n}{d}\oplus \bone}}\ar[d]^{T^\geq}_-{\cong} \\
 \ax{p-1}{\Pp^n} \times
 \ax{p}{\Pp^n} \times 
\ax{p}{\Pp^n} \ar[r]_-{\Psi_p} 
& \ax{p-1}{\Pp^n} \times \ax{p}{\Pp^{n}}
} 
$$
and sends
$
\Psi_p : \bigl(a[\Pp^{p-1}], b[\Pp^p], c[\Pp^p] \bigr) \longmapsto
\bigl( a[\Pp^{p-1}]+c\cdot c_{1}(\vb{\Pp^n}{d})\cap [\Pp^p] ,\ (b+c)[\Pp^p] \bigl)$.
Since \,
$c_{1}(\vb{\Pp^n}{d})\cap [\Pp^p]\, =\, d[\Pp^{p-1}]$\,
we then identify \,$\Psi_p$ with $\ (a, b, c) \mapsto (a+c\cdot d,\ b+c)$.

One can  associate variables $t_0 , \ t_1$ to the 
generators of $\ax{p}{\pb{\vb{\Pp^n}{d}\oplus \bone}} \cong \Zplus \oplus \Zplus$, and identify
$(\alpha_0 , \alpha_1) \in \ax{p}{\pb{\vb{\Pp^n}{d}\oplus \bone}}$ with
$t_0^{\alpha_0} t_1^{\alpha_1}$. It follows from Theorem~\ref{thm:main1}
and Remark~\ref{rem:power} that \ecs{p}{\pb{\vb{\Pp^n}{d}\oplus \bone}} 
can be written as
\begin{equation*}
\begin{split}
E_p & ( \Pp ( {\mathcal{O}}_{{\Pp}^n}(d) \oplus {\bone} ) ) = \\
&= \sum_{\alpha_{0},\alpha_{1}}
        \ecs{p}{\pb{\vb{\Pp^n}{d}\oplus \bone}} (\alpha_0 , \alpha_1) 
     t^{\alpha_0}_0 t^{\alpha_1}_1  \\
&= \sum_{\alpha_{0},\alpha_{1}}
   {\Psi_{p}}_{\#}\lp \ecs{p-1}{\pb{\vb{\Pp^n}{d}}} \odot \ecs{p}{\pb{\vb{\Pp^n}{d}}} 
   \odot \ecs{p}{\Pp^n} \rp (\alpha_0,\alpha_1) t^{\alpha_0}_0 t^{\alpha_1}_1  \\
&= \sum_{\alpha_{0},\alpha_{1}}\lp \sum_{(a,b,c)\in\Psi_{p}^{-1}(\alpha_0,\alpha_1)} 
    \ecs{p-1}{\Pp^n}(a)\cdot\ecs{p}{\Pp^n}(b)\cdot\ecs{p}{\Pp^n}(c)\rp
     t^{\alpha_0}_0 t^{\alpha_1}_1  \\
&= \sum_{\alpha_{0},\alpha_{1}}\lp 
    \sum_{\substack{a+cd=\alpha_0 \\ b+c=\alpha_1}} 
    \ecs{p-1}{\Pp^n}(a)\cdot\ecs{p}{\Pp^n}(b)\cdot\ecs{p}{\Pp^n}(c)\rp
    t^{a+cd}_0 t^{b+c}_1  \\
&= \lp\sum_{a\geq0} \ecs{p-1}{\Pp^n}(a) \cdot t_0^a\rp
    \lp\sum_{b\geq0} \ecs{p}{\Pp^n}(b)\cdot t_1^b\rp
    \lp\sum_{c\geq0}\ecs{p}{\Pp^n}(c)\cdot(t_0^{d}t_1)^c\rp.  
\end{split}
\end{equation*}
In Lawson-Yau~\cite{law&yau-hosy} it was shown that
\begin{equation*}
\sum_{k\geq0} \ecs{p}{\Pp^n}(k)\,t^k \, = \, 
        \lp\frac{1}{1-t}\rp^{\lp\substack{n+1\\p+1}\rp}
\end{equation*}
\noindent 
and hence \ecs{p}{\pb{\vb{\Pp^n}{d}\oplus \bone}}\ is then written as:
\begin{equation}
\label{2.21}
\ecs{p}{\pb{\vb{\Pp^n}{d}\oplus \bone}} \, = \, 
 \lp\frac{1}{1-t_0}\rp^{\lp\substack{n+1\\p}\rp}
 \lp\frac{1}{1-t_1}\rp^{\lp\substack{n+1\\p+1}\rp}\
 \lp\frac{1}{1-t_0^d t_{1}}\rp^{\lp\substack{n+1\\p+1}\rp}
\end{equation}
\end{example}
\begin{subexamples}\nopagebreak 
\label{subexmp:misc}
\begin{enumerate} 
\item When $d=0$ one has $\pb{\vb{\Pp^n}{d}\oplus \bone} = \Pp^n \times \Pp^1$ and 
our computations recover the formula
\begin{equation*}
\ecs{p}{\Pp^n\times\Pp^1} \, = \, 
\lp \frac{1}{1-t_0}\rp^{\lp\substack{n+1\\p}\rp} \ 
\lp\frac{1}{1-t_1}\rp^{2\lp\substack{n+1\\p+1}\rp}
\end{equation*}
obtained in~\cite{law&yau-hosy}. Furthermore, our general formula in 
Theorem~\ref{thm:main1} gives an inductive process to compute
$\ecs{p}{\Pp^n\times\Pp^m}$ in general.
\item When $n=1$, the $\Pp^1$-bundle \pb{\vb{\Pp^n}{d}\oplus \bone}\ over $\Pp^1$ is
the Hirzebruch surface $\F_{d}$, and the formula
\begin{equation*}
\ecs{1}{\F_{d}} \, = \, 
\lp\frac{1}{1-t_0}\rp^{2}\, \lp\frac{1}{1-t_1}\rp\,
\lp\frac{1}{1-t_0^d t_{1}}\rp
\end{equation*}
recovers the one obtained in Elizondo~\cite{eli-tor}.
\item When $d=1$, the $\Pp^1$-bundle 
\pb{\vb{\Pp^n}{1}\oplus \bone}\ over $\Pp^n$
is just the blow-up ${\widetilde{\Pp}}^{n+1}$ of $\Pp^{n+1}$ at
a point, and the expression
\begin{equation*}
\ecs{p}{ \widetilde{\Pp}^{n+1} } \, = \, 
\lp \frac{1}{1-t_0}\rp^{\lp\substack{n+1\\p}\rp}
 \lp\frac{1}{1-t_1}\rp^{\lp\substack{n+1\\p+1}\rp}\
 \lp\frac{1}{1-t_{0}t_1}\rp^{\lp\substack{n+1\\p+1}\rp}
\end{equation*}
recovers the formula obtained in~\cite{eli-tor}.
\end{enumerate}
\end{subexamples}

\section{Chow quotients and Euler-Chow series}
\label{sec:cq}

In this section we consider a projective algebraic variety $X$ equipped  
with an algebraic action of $\C^*$.
In general, an algebraic action of a torus $\T = (\C^*)^k$ 
provides stratifications on $X$, and here
we introduce the following one, much in the spirit of Kapranov \cite{kap-cq}.
 
\begin{definition}
\label{defn:type}
Let $\T=(\C^*)^k$ act algebraically on $X$
and let $\alpha\in \ax{p}{X}$ be a fixed
element with $0\leq p \leq k$. We say that a $p$-dimensional 
orbit $\T\cdot x$ is of {\em type} $\alpha$ if its closure
$\overline{\T\cdot x}$ lies in the component $\cvpd{p}{\alpha}{X}$,
when viewed as a $p$-dimensional effective cycle. One can stratify $X$
according to the orbit type of its elements, introducing the {\em Chow
stratification} of $X$.
\end{definition}

Assume that $\T$ has orbits of maximal dimension $k$.
It is easy to see that when $X$ is irreducible then there is a unique
maximal open stratum $X^o$ consisting of points whose orbits have maximal
dimension $k$ and are of the same type $\alpha_0$, for some $\alpha_0
\in \ax{k}{X}$. In particular, there is an embedding $X^o / \T
\hookrightarrow \cvpd{k}{\alpha_0}{X}$ and following
\cite{kap&stu&zel-quot} we introduce the following notion.

\begin{definition}
\label{defn:cq}
The {\em Chow quotient} $X\cq \T$ is the weak normalization of the
closure of $X^o/\T$ in $\cvpd{k}{\alpha_o}{X}$. 
\end{definition}
\begin{remark}
\label{rem:proper}
When $X^o/\T$ is proper, the Chow quotient gives a closed embedding of $X^o /\T$
into the appropriate Chow variety.
\end{remark}

We combine the use of Chow quotients with the notion of {\em trace
maps} introduced in Friedlander and Lawson \cite[\S 7.1]{fri&law-cocyc}, 
which is described as follows.
A Chow variety $\cvpd{k}{\alpha}{X}$, being projective, 
has its own Chow monoid
$\cvp{p}{\cvpd{k}{\alpha}{X}}$. 
The trace map is a continuous  monoid morphism
\begin{equation}
\label{eqn:tr-law}
tr: \cvp{p}{\cvpd{k}{\alpha}{X}} \rightarrow \cvp{p+k}{X}
\end{equation}
that, roughly speaking, associates to an irreducible $p$-dimensional  
subvariety $W\subset \cvpd{k}{\alpha}{X}$ its ``total fundamental cycle'',
which is a $(p+k)$-cycle in $X$. 

Combining the two constructions above,  
one associates to such a $\T$-action on an irreducible variety $X$,  
a trace map
\begin{equation}
\label{eqn:trace}
t_p : \cvp{p}{X\cq \T} \rightarrow
           \cvp{p+k}{X}.
\end{equation}
This map, in turn, induces a monoid morphism
\begin{equation}
\label{eqn:trace0}
\varphi_p : \ax{p}{X\cq \T} \rightarrow
           \ax{p+k}{X},
\end{equation}
in the level of $\pi_0$. In other words,
$t_p ( \cvpd{p}{\alpha}{X\cq \T} ) \subset 
\cvpd{p+k}{\varphi_p(\alpha)}{X}.$
Explicit examples of such monoid morphisms are given in
Proposition \ref{prop:trace1} and Corollary \ref{cor:tr-flag} below.

\begin{proposition}
\label{prop:inj}
Let $\T\cong (\C^*)^k$ act on an irreducible projective variety
$X$ with generic orbits of maximal dimension $k$. If
$\ \dim{(X\cq \T - X^o/\T )} < p \leq \dim{X} - k \ $ then
the trace map $t_p : \cvp{p}{X\cq \T} \rightarrow \cvp{p+k}{X}$
is injective.
\end{proposition}
\begin{proof}
Let $W \subset X\cq \T$ be an irreducible variety of dimension
$p$. The  hypothesis on dimensions implies that  $W^o \equdef W \cap
X^o/\T$ is an open dense subvariety of $W$. In order to define the
trace $t_p(W)$ one considers the cycle $Z_W \subset W \times X$, 
defined as the ``closure of the cycle'' 
$$Z_W^o = \{ (w,c_w) \ | \ w\in W^o \text{ and } 
  c_w = \text{ cycle whose Chow point is } w \};$$ 
cf. \cite[\S 7.1]{fri&law-cocyc}. 
By definition, 
$t_p(W) = pr_{2*}(Z_W)$, 
where $pr_2$ is projection onto the second factor of $W\times X$.

Note that, since the fibers of the projection $Z_W^o \to W^o$ are
irreducible and $W^o$ is irreducible, then $Z_W$ is irreducible of
dimension $p+k \leq \dim{X}$. Furthermore, if $(a,b), (a',b)
\in Z_W^o$ then $a=a'$, for if $a,a'\in W^o \subset X^o/\T$ 
then the orbits in $X$ represented by $a$ and $a'$
have a common point $b$, hence $a=a'$. It follows that 
$pr_2 $ maps $Z_W$ birationally onto its image, and thus
$t_p(W) \equdef pr_{2*} ( Z_W )$ is an irreducible cycle 
with multiplicity $1$.

Suppose that $t_p(W) = t_p(W')$, where both $W$ and $W'$ are
irreducible. An element $w\in W^o$ corresponds to an orbit
$\T\cdot x \subset pr_2(Z_W) = pr_2(Z_{W'})$, with $x\in X^o$.
Therefore $w\in {W'}^o$, and hence $W=W'$. Since that Chow
monoids are freely generated by the irreducible subvarieties, one
concludes that the trace map $t_p$ is injective.
\end{proof}

We now consider the situation where
$X$ is an irreducible projective variety, 
on which $\C^*$ acts in such a fashion that the fixed point 
set $X^{\C^*}$ has only two connected components $X_1$ and $X_2$.  
Let $i_1 : X_1 \hookrightarrow X$ and
$i_2 : X_2 \hookrightarrow X$ denote the inclusion maps, and let 
$t_{p-1} : \cvp{p-1}{X\cq \C^*} \rightarrow \cvp{p}{X}$ be the trace
morphism (\ref{eqn:trace0}). These maps induce a monoid morphism
\begin{align}
\label{eqn:fund0}
\Psi_p : \ax{p-1}{X\cq \C^*}\times  \ax{p}{X_1} \times \ax{p}{X_2}
&\rightarrow \ax{p}{X} \\
(\alpha,\beta,\gamma) & \longmapsto \varphi_{p-1}(\alpha) + {i_1}_*\beta +
 {i_2}_* \gamma, \notag
\end{align}
which yields our next result.

\begin{theorem}
\label{thm:cq}
Let $X$ be an smooth projective  variety 
on which $\C^*$ acts algebraically.
If $X^{\C^*}$ is the union of two connected components $X_1$ and $ X_2$,
then for each $0\leq p\leq \dim{X}$ one has
\begin{equation}
\ecs{p}{X} = {\Psi_p}_\sharp \lp \ecs{p-1}{X\cq \C^*}\odot
                            \ecs{p}{X_1} \odot \ecs{p}{X_2} \rp.
\end{equation}
\end{theorem}
\begin{proof}
The constructions above give a continuous monoid morphism
\begin{align}
\label{eqn:fund1}
\phi_p : \cvp{p-1}{X\cq \C^*}\times  \cvp{p}{X_1} \times \cvp{p}{X_2}
&\rightarrow \cvp{p}{X} \\
(a,b,c) & \longmapsto t_{p-1}(a) + {i_1}_*b +
 {i_2}_* c, \notag
\end{align}
which induces $\Psi_p$ in the level of $\pi_0$; cf. (\ref{eqn:fund0}). 

It is clear that the image of $\phi_p$ lies
in the fixed point set $\cvp{p}{X}^{\C^*}$.
We show that $\phi_p$ is injective and that it surjects onto
$\cvp{p}{X}^{\C^*}$. 

Suppose that $\phi_p(a,b,c)=\phi_p(a',b',c')$. By ``disjointness
of support'' one sees immediately that 
\begin{equation}
\label{eqn:varphi}
b=b',\quad c=c'\quad \text{and}\quad t_{p-1}(a) = t_{p-1}(b).
\end{equation}
 
It follows from the decomposition of $X$ described in Bialynicki-Birula
\cite[\S 4]{bb-red} that $X^o = X - X^{\C^*}, \ $ and that
$X^o/\C^*$ is proper. Therefore $X\cq \C^* = X^o / \C^*,\ $ and
Proposition \ref{prop:inj} implies that the trace map $t_{p-1}$ is
injective, which together with (\ref{eqn:varphi}) shows that
$\phi_p$ is injective. 

Given an irreducible, $p$-dimensional, and $\C^*$-invariant 
subvariety $V\subset X$, then either $V \subset X^{\C^*}$, in
which case it clearly lies in the image of $\phi_p$, or $V^o =
V\cap X^o \neq \emptyset$. In the latter case, define 
$W = \overline{V^o/\C^*} \subset X\cq \C^*$. A simple
inspection shows that $t_p(W)  = V$. This suffices to show that
$\phi_p$ surjects onto $\cvp{p}{X}^{\C^*}$ and, using 
Lemma \ref{lem:tech}, one concludes that $\phi_p$ is a homeomorphism
onto its image.

The rest of the proof uses the same arguments as in the last
paragraph in the proof of Theorem \ref{thm:main1}.
\end{proof}

\begin{remark}
\label{rem:proj-clo}
One can directly prove that the Chow quotient $\Pp(E_1\oplus E_2)\cq
\C^*$ is isomorphic to $\Pp(E_1)\times_W \Pp(E_2)$, whenever $E_1$ and
$E_2$ are algebraic vector bundles over a variety $W$.
In this case, Theorem \ref{thm:main1} becomes a consequence of the
result above whenever $W$ is smooth.
\end{remark}

Next, we describe some examples of Chow quotients, trace maps 
and resulting computations of Euler-Chow series.

\subsection{Examples} \hfill
\label{subsec:exmp}

The linear action of $\C^*$ on the last coordinate of $\C^{n+1}$
induces an algebraic $\C^*$ action on the Grassmannian $\gr{d}{n}$
of $d$-planes in $\Pp^n$. More generally, one obtains an algebraic
action on all partial flag varieties $\fl{d_1,\ldots,d_r}{n}$ 
of nested linear subspaces $D_1\subset \cdots \subset D_r$ of $\Pp^n$, 
satisfying $\dim{D_i} = d_i$.

We first describe the orbit structure and the Chow quotient
of both  $\gr{d}{n}$ and $\fl{0,1}{n}$, under this $\C^*$
action.
\bigskip

\subsubsection{The Grassmannian case}

We adopt the convention that $\gr{d}{n}= \emptyset$, whenever
$d<0$, and that $\chi(\emptyset) =0.$
Fix $\py = [0: \cdots : 0:1]\in \Pp^n$, and let 
$L\in \gr{d}{n}$ be fixed under the above action. Then
the corresponding $d$-dimensional subspace $L\subset \Pp^n$ can
be of two types:
\begin{enumerate}
\item Either $L$ is contained in 
$\Pp^{n-1} = \{ [x_0: \cdots : x_n] \in \Pp^n \ | \ x_n = 0 \}$,
\item Or $L$ has the form $\py \# \ell$, where $\ell$ is a
$(d-1)$-dimensional subspace of $\Pp^{n-1}$, and $\py\# \ell$ denotes
the ruled join of $\py$ and $\ell$ in $\Pp^{n}$.
\end{enumerate}

In other words, the fixed point set $\gr{d}{n}^{\C^*}$ has two
connected components which are naturally isomorphic to 
$\gr{d}{n-1}$ and $\gr{d-1}{n-1}$, and whose
inclusions in $\gr{d}{n}$ are denoted by 
$i: \gr{d}{n-1} \hookrightarrow \gr{d}{n}$ and $j: \gr{d-1}{n-1}
\hookrightarrow \gr{d}{n}$, respectively. In particular,
this $\C^*$ action on $\gr{d}{n}$  satisfies the hypothesis of
Theorem \ref{thm:cq}. Furthermore, all points in the generic locus
$\gr{d}{n}^o = \gr{d}{n} - \gr{d}{n}^{\C^*}$ 
have the same orbit type.

\begin{proposition}
The Chow quotient $\gr{d}{n}\cq \C^*$ is naturally isomorphic to
the flag variety $\fl{d-1,d}{n-1}$.
\end{proposition}
\begin{proof}
Let $\pi : \Pp^n - \{ p_\infty \} \rightarrow \Pp^{n-1}$
denote the projection onto the hyperplane $\Pp^{n-1}$. 
Standard arguments show that
\begin{align}
q: \gr{d}{n}^o & \rightarrow \fl{d-1,d}{n-1} \\
L & \mapsto ( L\cap \Pp^{n-1}, \pi(L) ) \notag
\end{align}
is a regular, surjective map. This map descends to the quotient
$\gr{d}{n}^o /\C^*$ and produces a closed inclusion
$\gr{d}{n}\cq \C^* \equiv \gr{d}{n}^o / \C^* \hookrightarrow
\fl{d-1,d}{n-1}$; cf. Proposition \ref{prop:inj} and Remark
\ref{rem:proper}. This is then easily seen to be an isomorphism.
\end{proof}
\begin{corollary}
\label{cor:grass}
The $p$th Euler-Chow series of the Grassmannian $\gr{d}{n}$ is given by 
$$\ecs{p}{\gr{d}{n}} = \Psi_{p\sharp} \lp
\ecs{p-1}{\fl{d-1,d}{n-1}} \odot
\ecs{p}{\gr{d}{n-1}} \odot
\ecs{p}{\gr{d-1}{n-1}} \rp,
$$
where $\Psi_p$ is given by {\it (\ref{eqn:fund0})}.
\end{corollary}

We proceed to explicitly describe $\Psi_p$ in this example.
We follow closely the projective notation for the Schubert cycles in
Grassmanians and flag varieties, as described in Fulton
\cite[\S 14.7]{ful-inter}. However, since we are dealing with various
projective spaces, we add an upperscript in the notation to
denote the dimension of the ambient projective space.

First, fix a complete flag
$L_0 \subset L_1 \subset \cdots \subset L_n = \Pp^n$ of linear
subspaces, and associate to a sequence of length $d$  
$\mba : 0\leq a_0 < a_1 < \cdots < a_d\leq n$, the Schubert variety
$\sch{\mba}{n} = \sch{a_0,\ldots,a_d}{n}\subset \gr{d}{n}$ 
defined by
$$\sch{\mba}{n} = 
\{ l \in \gr{d}{n} \ | \ \dim{(l \cap L_i)} \geq a_i, i=0,\ldots d
\}.$$

Now, let $\{ \mba^1\subset \ldots \subset \mba^r \}$ be a nested
collection of sequences, 
such that the length of $\mba^i$ is $d_i$, and define the 
Schubert variety
$\sch{\mba^1;\ldots;\mba^r}{n} \subset \fl{d_1,\ldots,d_r}{n}$
by
$$
\sch{\mba^1;\ldots;\mba^r}{n} = 
\{ (l_1,\ldots,l_r) \in \fl{d_1,\ldots,d_r}{n} \ | \ 
l_i \in \sch{\mba^i}{n} \}.$$

\begin{facts}\hfill
\begin{enumerate}
\item The Schubert variety $\sch{\mba^1;\ldots;\mba^r}{n}$ 
is an irreducible subvariety of $\fl{d_1,\ldots,d_r}{n}$ 
of dimension 
$d(\mba^1;\ldots;\mba^r) = \sum_{i=1}^r {\sum_{j=0}^{'d_i}}(a_j^i -j)$.
In this sum, any term $(a_j^i - j)$ is omitted if $a_j^i$ appears in
the previous $\mba^{i-1}$.
\item The surjection $\ax{p}{\fl{d_1,\ldots,d_r}{n}} \rightarrow
\algp{p}{\fl{d_1,\ldots,d_r}{n}}$ is an isomorphism. Therefore, if
$\scho{\mba^1;\ldots;\mba^r}{n}$ denotes the class of 
$\sch{\mba^1;\ldots;\mba^r}{n}$ in
$\ax{*}{\fl{d_1,\ldots,d_r}{n}}$, then the monoid 
$\ax{p}{\fl{d_1,\ldots,d_r}{n}}$ is
freely generated by the collection
$$\{ \scho{\mba^1;\ldots;\mba^r}{n} \ | \ d(\mba^1;\ldots;\mba^r) = p \}.$$ 
\end{enumerate}
\end{facts}

\begin{remark}
Let $\mba$ be the sequence $\ 0<1< \cdots < d-1< d+1.$
Then, it is easy to see that the class of a $1$-dimensional orbit
$[\overline{\T\cdot x}] \in \ax{1}{\gr{d}{n}}$, 
is precisely $\scho{\mba}{n}$.
\end{remark}

Consider the universal bundles $U_{d}$ and $U_{d+1}$ over
the flag variety $\fl{d-1,d}{n-1}$, of ranks $d$ and $d+1$, respectively.
Then, denote the quotient line bundle $U_{d+1}/U_d$ by $L$, and let 
$\rho : \Pp(L\oplus \bone) \rightarrow \fl{d-1,d}{n-1}$ be the bundle
projection. The following result describes the trace map
$t_p : \cvp{p}{\fl{d-1,d}{n-1}} \rightarrow 
     \cvp{p+1}{\gr{d}{n}}$, 
in a very explicit and geometric manner.
\vspace{.36in}

\begin{proposition}\hfill
\label{prop:trace1}

\begin{enumerate}
\item The blow-up of $\gr{d}{n}$ along the fixed point set
$\gr{d}{n}^{\C^*}= \gr{d-1}{n-1} \amalg \gr{d}{n-1}$ is the variety $ \Pp(L \oplus
\bone)$.
\item Let $b: \Pp(L\oplus \bone) \rightarrow \gr{d}{n}$ denote the
blow-up map. Then $\phi_p : \cvp{p}{\fl{d-1,d}{n-1}} \rightarrow
\cvp{p+1}{\gr{d}{n}}$ is the composition $b_* \circ \rho^*,$ where
$b_*$ denote the proper push-forward under $b$, and $\rho^*$ is the flat
pull-back. 
\end{enumerate}
\end{proposition}
\begin{proof}
Left to the reader.
\end{proof}

The induced map 
$\varphi_p : \ax{p}{\fl{d-1,d}{n-1}} \rightarrow \ax{p+1}{\gr{d}{n}}$
can be computed either from the proposition or by a direct
argument.

\begin{corollary} 
\label{cor:tr-flag}
Given a sequence 
$0\leq a_1 <\cdots <a_d \leq n-1$ and $0\leq j \leq d$, one has:
\begin{equation}
\label{eqn:schubert}
\varphi_p ( \scho{a_0,\ldots,\widehat{a_j},\ldots,a_d;
        a_0,\ldots,a_d}{n-1} ) = 
        \scho{a_0,\ldots,a_{j-1},a_j+1,\ldots,a_d+1}{n}.
\end{equation}
\end{corollary}
\begin{proof}
Exercise for the reader.
\end{proof}

\subsubsection{The flag varieties $\fl{0,1}{n}$}

The orbit structure of the $\C^*$ action on $\fl{0,1}{n}$ is
described in the following statement. Here we denote by  
$Q_{d,n}$ and $S_{d,n}$ the  universal quotient bundle and the
tautological bundle over $\gr{d}{n}$, respectively.

\begin{proposition}\hfill
\label{prop:flag-orbit}
\begin{enumerate}
\item The fixed point set $\fl{0,1}{n}^{\C^*}$ has three connected
components $F_1, F_2$ and $F_3$, respectively isomorphic to $\fl{0,1}{n-1}$,
$\Pp^{n-1}$ and $\Pp^{n-1}$.
\item There are two $\C^*$-invariant subvarieties $W_1$ and $W_2$
of $\fl{0,1}{n}$ which are equivariantly isomorphic to the projective bundles
$\Pp(Q_{0,n-1} \oplus \bone)$ and $\Pp( S_{0,n-1} \oplus \bone)$
over $\Pp^{n-1}$,  respectively. 
Furthermore, $W_1$ and $W_2$ are Schubert
cycles for an appropriate choice of coordinates, such that their
Schubert classes are given by
$[W_1]= \scho{n-1;n-1,n}{n}$ and $[W_2]= \scho{n;0,n}{n}$.
\item The component $F_1$ is precisely $\Pp(Q_{0,n-1})
\subset \Pp(Q_{0,n-1}\oplus \bone)$. 
The strata $W_1$ and $W_2$ intersect at 
$F_2$, where $F_2 \equiv \Pp( \bone ) \subset \Pp(Q_{0,n-1}\oplus \bone)$
and $F_2 \equiv \Pp( S_{0,n-1} )\subset \Pp(S_{0,n-1}\oplus \bone)$.
The component $F_3$ is identified with 
$\Pp(\bone)\subset \Pp(S_{0,n-1}\oplus \bone)$.
\end{enumerate}
\end{proposition}
\begin{proof}
The first assertion is clear.

To prove the second assertion, one needs only to consider the
$\C^*$-equivariant projections 
$\pi_1 : \fl{0,1}{n}= \Pp ( Q_{0,n} )   \rightarrow \gr{0}{n} \equiv \Pp^n$
and 
$\pi_2 : \fl{0,1}{n}=\Pp( \check{S}_{1,n} ) \rightarrow \gr{1}{n}$.
Then $W_1$ is defined to be $\pi_1^{-1}( \Pp^{n-1} ) \subset
\fl{0,1}{n}$ and $W_2 = \pi_2^{-1} ( \gr{0}{n-1} ) \subset
\fl{0,1}{n}$, where $\gr{0}{n-1}$ is one of the components of 
$\gr{1}{n}^{\C^*}$ and $\Pp^{n-1} = \{ x_n = 0 \}$. 
The rest of the proof is an easy consequence of
this description
\end{proof}

Applying \cite[Lemma 4.1]{bb-red} to our situation, one concludes 
that the closure of a  $1$-dimensional orbit $\C^*\cdot x$ in $\fl{0,1}{n}$
intersects precisely $2$ connected components of the fixed point set.
Furthermore, an application of \cite[Theorems 4.1 and 4.3]{bb-red}
implies that these two components completely determine the {\em type}
of the orbit; cf. Definition~\ref{defn:type}.

There are 3 types of $1$-dimensional orbits:
\begin{enumerate}
\item {\bf Type 1:} The $1$-dimensional orbits contained in
the stratum $W_1$. These are the orbits of the points in $W_1 - \{ F_1
\cup F_2 \}$.
\item {\bf Type 2:} The $1$-dimensional orbits contained in
the stratum $W_2$. These are the orbits of the points in $W_2- \{ F_2
\cup F_3 \}$.
\item {\bf Generic type:} These are the orbits of points in the
{\em generic locus} $\fl{0,1}{n}^o = \fl{0,1}{n} - \{ W_1 \cup W_2\}$
of the action.
\end{enumerate}
\begin{proposition} \hfill

\begin{enumerate}
\item The Chow quotient $W_1 \cq \C^*$ is naturally isomorphic to
$F_1$.
\item The Chow quotient $W_2 \cq \C^*$ is naturally isomorphic to
$F_3$.
\item The Chow quotient $\fl{0,1}{n}\cq \C^* $ is naturally isomorphic
to $\Pp^{n-1}[2]$, the blow-up of $\Pp^{n-1}\times \Pp^{n-1}$ along
the  diagonal.
\end{enumerate}
\end{proposition}
\begin{proof}
The two first assertions are simple and follow from 
Remark \ref{rem:proj-clo}.
In order to prove the last assertion, 
we first identify $\Pp^{n-1}[2]$ with the fibered product 
$$\fl{0,1}{n-1}\times_{\gr{1}{n-1}}\fl{0,1}{n-1} = \{ (l,l',L)\in
\Pp^{n-1} \times \Pp^{n-1} \times \gr{1}{n-1} \ | \ l,l'\subset L
\}.$$ Using the notation of Proposition \ref{prop:flag-orbit} we denote
$\fl{0,1}{n}^o = \fl{0,1}{n}- \{ W_1\cup W_2 \}$ and 
define an algebraic map
\begin{align}
 \tau : \fl{0,1}{n}^o & \rightarrow \Pp^{n-1}[2] \\
(l, L) &\mapsto (l\cap \Pp^{n-1}, \pi(l), \pi(L)). \notag
\end{align}
Notice that this map factors through $\fl{0,1}{n}^o/\C^*$,
defining
$\tau' : \fl{0,1}{n}^o/\C^* \rightarrow \Pp^{n-1}[2]$, and that it actually
defines an isomorphism 
$\tau' : \fl{0,1}{n}^o/\C^* \rightarrow \Pp^{n-1}[2]- D$, where 
$D= \{ (l,l,L) \ | \ l\in L \}$ is the
exceptional divisor of $\Pp^{n-1}[2]$. Let 
$\delta' : \Pp^{n-1}[2] - D \rightarrow \fl{0,1}{n}^o/\C^*$ denote the
inverse of this isomorphism. Since $\Pp^{n-1}[2]$ is smooth and
the Chow quotient $\fl{0,1}{n}\cq \C^*$ is irreducible and weakly normal, 
one concludes that $\delta'$ has a unique extension to a 
surjective map 
\begin{equation}
\label{eqn:delta}
\delta : \Pp^{n-1}[2]  \to \fl{0,1}{n}\cq \C^*.
\end{equation}

In order to prove that $\delta$ is an isomorphism, one just needs to 
prove that it is injective, due to the properties of birational maps
between weakly normal spaces; cf. Koll\'ar \cite{kol-rat}. The injectivity will
follow from the following description of $\delta$.

Let $\Gamma \subset \Pp^{n-1}[2] \times \fl{0,1}{n}$ be defined as the
closure of 
\begin{equation}
\label{eqn:gamma}
\Gamma^o = \left\{ (a,b,C)\times (l,L) \ \left| 
\ (a,b,C) \in \Pp^{n-1}[2] - D, \
(l,L)\in \fl{0,1}{n}, \
L \cap \Pp^{n-1} = a,  \
\pi(l) = b  
\right.  \right\} .
\end{equation}
We observe that the projection $q : \Gamma \to \Pp^{n-1}[2]$ is
equidimensional. Indeed, it is easy to see that
if $(l,l',L)\in \Pp^{n-1}[2] -D$ then $q^{-1}(l,l',L) = \supp{\delta
(l,l',L)}$. In other words,  $q^{-1}(l,l',L)$ is the closure of the orbit
of $( a , l \# a )$, where $a$ is any point in $l'\# \py - \{ l',\py \}$.
Now, given $(l,l,L)\in D$, 
one can think of $(l,L)$ as an element in
$F_1 \subset \fl{0,1}{n}$, and of $l$ as an element in $F_3 \subset
\fl{0,1}{n}$. It follows from the first two assertions of this
proposition that there are $1$-dimensional orbits
$\psi_1(l,L) \subset W_1$ and $\psi_2(l)\subset W_2$ uniquely determined
by $(l,l,L)$. It is easy to see that $q^{-1}(l,l,L) =
\psi_1(l,L)\cup \psi_2(l)$.

It follows from \cite{fri&law-cocyc} that the equidimensional 
projection $q: \Gamma \to \Pp^{n-1}[2]$ induces an
algebraic map $\ q_\Gamma : \Pp^{n-1}[2] \to \cvpd{1}{\alpha}{\fl{0,1}{n}}$.
This map concides with $\delta$ on
$\Pp^{n-1}[2] -D$ and hence,  by uniqueness, one concludes that $\delta =
q_\Gamma$. In this case, $\delta(l,l,L) = \psi_1(l,L)+\psi_2(l)$, 
and hence  $\delta$ is injective. This concludes the proof.
\end{proof}

We now compute the Euler-Chow series for the flag variety
$\fl{0,1}{2}$, using the above information. In this case one has
$$F_1  \cong W_1\cq \C^*\cong F_2 \cong F_3\cong W_2\cq \C^* \cong
\Pp^1,$$
and 
$$\fl{0,1}{2} \cq \C^* \cong \Pp^1[2] = \Pp^1\times \Pp^1.$$

It follows from the isomorphisms
$\ax{1}{\fl{0,1}{2}}\cong \Z_+\cdot \scho{1;0,1}{2}
\oplus \Z_+ \cdot \scho{0;0,2}{2} \cong \Z_+ \oplus \Z_+$, that we may
denote the connected components of the Chow monoid $\cvp{1}{\fl{0,1}{2}}$ by 
$\cvpd{1}{(r,s)}{\fl{0,1}{2}}$, where $(r,s)\in \Z_+\oplus \Z_+$.

The Chow quotients described above induce inclusions:
\begin{enumerate}
\item $\Pp^1 \cong W_1\cq \C^* \subset \cvp{1}{\fl{0,1}{2}}$,
\item $\Pp^1 \cong W_2\cq \C^* \subset \cvp{1}{\fl{0,1}{2}}$, and
\item $ \Pp^1\times \Pp^1 \cong \fl{0,1}{2}\cq \C^* \subset 
        \cvp{1}{\fl{0,1}{2}}$,
\end{enumerate}
whose associated trace maps are given as follows.
\begin{lemma} \hfill
\label{lem:flag}

\begin{enumerate}
\item $\psi_1 : \cvpd{1}{d}{\Pp^1} \equiv \cvpd{1}{d}{W_1\cq \C^*}\rightarrow
\cvpd{2}{(d,0)}{\fl{0,1}{2}}$, sends $d\cdot \Pp^1 $ to
$d\cdot W_1$;
\item $ \psi_2 : \cvpd{1}{d}{\Pp^1} \equiv \cvpd{1}{d}{W_2\cq \C^*}\rightarrow
\cvpd{1}{(0,d)}{\fl{0,1}{2}}$, sends $d\cdot \Pp^1 $ to
$d\cdot W_2$;
\item Let 
$ t_1 : \cvpd{1}{(r,s)}{\Pp^1\times \Pp^1} \rightarrow
\cvpd{2}{(r,s)}{\fl{0,1}{2}} $ 
be the trace map induced by the 
Chow quotient $\Pp^1\times \Pp^1 = \fl{0,1}{2}\cq \C^*$,  and let $Z$ be an 
irreducible subvariety of $\Pp^1\times \Pp^1$. If the $2$-cycle
$t_1(Z)$ contains either $W_1$ or $W_2$ as an irreducible component,
then $Z = \Delta =$  the diagonal of  $\Pp^1\times \Pp^1$, 
in which case $t_1(\Delta) =  W_1+W_2.$
\end{enumerate}
\end{lemma}
\begin{proof}
The first two assertions are obvious from the definitions.
In order to prove the last assertion, consider an irreducible
$1$-dimensional subvariety $Z\subset \Pp^1\times \Pp^1$, distinct from
the diagonal $\Delta$, and assume that $W_1 \subset \supp{t_1(Z)}$. Since
$\supp{t_1(Z)} = \cup_{w\in Z} \supp{\delta(w)}$ and $Z\cap \Delta$ is finite,
one concludes that
$W_1 \cap \lp  \cup_{w\in Z-\Delta} \supp{\delta(w)} \rp$ 
is open and dense in $W_1$. 
On the other hand, if $w\in Z- \Delta$ then 
$\supp{\delta(w)}$ is the closure of a generic orbit, and hence it only
intersects $W_1$ at $F_1$. This implies that 
$W_1 \cap \lp  \cup_{w\in Z-\Delta} \delta(w) \rp \subset F_1$ is open in
$W_1$, which is a contradiction. Similar arguments are used in the case 
$W_2 \subset \supp{t_1(Z)}.$ One then concludes that $Z=\Delta$, and the
rest of the proof follows easily.
\end{proof}
With the above data, we can prove the following:
\begin{theorem}
\label{thm:f012div}
Let $x$ and $y$ be variables associated to the Schubert classes
$\scho{1;1,2}{2}$ and $\scho{2;0,2}{2}$, respectively.
Then the $2$nd Euler-Chow series of  $\fl{0,1}{2}$
is given by the generating function
\begin{equation}
\ecs{2}{\fl{0,1}{2}} = \frac{1-xy}{(1-x)^3(1-y)^3}.
\end{equation}
In other words,
$$
\ecs{2}{\fl{0,1}{2}}(r,s) \equdef \chi \lp \cvpd{2}{(r,s)}{\fl{0,1}{2}}\rp= 
\frac{1}{2}(r+1)(s+1)(r+s+2).
$$
\end{theorem}
\begin{proof}
Let
$\nu_{(r,s)}\ : \ \fl{0,1}{2}\ \to \
\Pp( \text{Sym}^r(\Lambda^1 \C^3)\otimes
     \text{Sym}^s(\Lambda^2 \C^3) ),
$
denote the composition
\begin{align*}
\fl{0,1}{2} & \hookrightarrow \gr{0}{2}\times \gr{1}{2} 
\xrightarrow{ p_1\times p_2 } \Pp( \Lambda^1\C^3 ) \times 
\Pp(\Lambda^2 \C^3) \\
& \xrightarrow{ \nu_r \times \nu_s } 
\Pp( \text{Sym}^r(\Lambda^1 \C^3) ) \times
\Pp( \text{Sym}^s(\Lambda^2 \C^3) ) \xrightarrow{ s }
\Pp( \text{Sym}^r(\Lambda^1 \C^3)\otimes
     \text{Sym}^s(\Lambda^2 \C^3) ),
\end{align*}
where the first map is the canonical inclusion, the second one is a
product of  Pl\"ucker embeddings, the third map is
the product of the appropriate Veronese maps (embeddings when
$r>0$ and $s>0$) and the last 
one is a Segre embedding. Let $\calo_{F}(r,s)$ denote the
pull-back $\nu_{(r,s)}^*(\calo(1))$ of the hyperplane bundle.
Then, it is easy to see that the Chow variety
$ \cvpd{2}{(r,s)}{\fl{0,1}{2}} $ is precisely the linear system
$ \Pp\left( H^0 ( \fl{0,1}{2}, \calo_{F}(r,s) ) \right)$. 
It follows from the Borel-Weil theorem that if $\C^3 $ denotes 
the canonical representation of $\text{GL(3,\C )}$, then
$ H^0 ( \fl{0,1}{2}, \calo_{F}(r,s) ) $
is the irreducible $\text{GL}(3,\C)$-module 
$\mathbb{S}_\lambda(\C^3)$ of highest weight $\lambda = ( r+s, s, 0)$.
See \cite[\S 9.3]{ful-tab} for complete details. 
Therefore, well-known formulas \cite[Theorem 6.3]{ful&har-rep} for the dimension of the Schur module
$\mathbb{S}_\lambda \C^3$ gives
$$
\chi \left( \cvpd{2}{(r,s)}{\fl{0,1}{2}} \right)
= \dim \mathbb{S}_\lambda( \C^3 ) =
\frac{1}{2}(r+1)(s+1)(r+s+2)
. 
$$
This concludes the proof.
\end{proof}
\begin{remark}
We thank the referee for suggesting this proof, and for the 
observation that all divisorial Euler-Chow functions for arbitrary
flag varieties should be given in terms of dimensions of irreducible
representations of $\text{GL}(n,\C)$. The various resulting
generating  functions  then should be given by well-known 
formulas; cf. \cite{mcd-poly}.
\end{remark}

In order to illustrate a more topological approach, using
techniques which can extended to more general situations,
we present our original proof below. 

{\small{ 
\begin{proof}
We use the notation of \ref{prop:flag-orbit}. 
Denote $F^0_{r,s}= \cvpd{2}{(r,s)}{\fl{0,1}{2}}$, and let 
$F^1_{r,s} \subset F^0_{r,s} $ be the (Zariski) closed subset
consisting of those effective divisors which contain either $W_1$ or
$W_2$ in their support. Note that 
\begin{equation} 
\label{eqn:union}
F^1_{r,s} = \lp W_1+ F^0_{r-1,s} \rp \cup \lp W_2 + F^0_{r,s-1} \rp
\end{equation}
and that
\begin{equation}
\label{eqn:intersection}
\lp W_1+ F^0_{r-1,s} \rp \cap \lp W_2 +
F^0_{r,s-1}\rp = (W_1+W_2)+ F^0_{r-1,s-1}.
\end{equation}

It follows from Lemma  \ref{lem:flag}(3) that the commutative diagram
\begin{equation}
\label{eqn:flag}
\begin{CD}
\Delta + \cvpd{1}{(r-1,s-1)}{\Pp^1\times \Pp^1} @>>>
\cvpd{1}{(r,s)}{\Pp^1\times \Pp^1} \\
@VVV @VVV \\
F^1_{r,s} @>>> F^0_{r,s}
\end{CD}
\end{equation}
is a fiber square whose horizontal arrows are closed inclusions, 
and such that the map of pairs 
\begin{equation}
\label{eqn:relative}
\left( \cvpd{1}{(r,s)}{\Pp^1\times \Pp^1}, \ 
  \Delta+\cvpd{1}{(r-1,s-2)}{\Pp^1\times \Pp^1} \right) 
\longmapsto 
 (F^0_{r,s},\ F^1_{r,s})
\end{equation}
is a relative homeomorphism. 
In order to simplify notation, let us write
\begin{equation}
a^0_{r,s} = \chi( F^0_{r,s}),\quad  \ a^1_{r,s}= \chi( F^1_{r,s}) 
\quad \
\text{and} 
\quad 
b_{r,s} = \chi\left( \cvpd{1}{r,s}{\Pp^1\times \Pp^1} \right) .
\end{equation}
It follows from
(\ref{eqn:flag}) and (\ref{eqn:relative}), and the additivity properties
of the Euler characteristic, that
\begin{equation}
\label{eqn:euler}
a^0_{r,s} =  a^1_{r,s} + b_{r,s} - \ b_{r-1,s-1}.
\end{equation}
Furthermore, (\ref{eqn:union}) and (\ref{eqn:intersection}) imply that
\begin{equation}
\label{eqn:a1}
a^1_{r,s} = a^0_{r-1,s}+a^0_{r,s-1}-a^0_{r-1,s-1}.
\end{equation}

With $x$ and $y$ as in the statement of the theorem, one has
\begin{align*}
E =& \ \ecs{2}{\fl{0,1}{2}}  
= \sum_{r\geq 0,\ s\geq 0} a^0_{r,s} x^ry^s \\
=& \sum_{r\geq 0,\ s\geq0} \lp
   a^1_{r,s} + b_{r,s} - b_{r-1,s-1} \rp x^ry^s \\
=& \sum_{r\geq 0,\ s\geq0} \lp
   a^0_{r-1,s} + a^0_{r,s-1} - a^0_{r-1,s-1}
   + b_{r,s} - b_{r-1,s-1} \rp x^ry^s \\
=& \sum_{r\geq 0,\ s\geq0} a^0_{r-1,s}   x^ry^s + 
   \sum_{r\geq 0,\ s\geq0} a^0_{r,s-1}   x^ry^s - 
   \sum_{r\geq 0,\ s\geq0} a^0_{r-1,s-1} x^ry^s    \\
&+ \sum_{r\geq 0,\ s\geq0} b_{r,s}       x^ry^s - 
   \sum_{r\geq 0,\ s\geq0} b_{r-1,s-1}   x^ry^s     \\
=& \sum_{r\geq 0,\ s\geq0} a^0_{r-1,s}   x^ry^s + 
   \lp \sum_{r\geq 0,\ s\geq0} a^0_{r,s-1}   x^ry^{s-1}\rp y - 
   \lp \sum_{r\geq 0,\ s\geq0} a^0_{r-1,s-1} x^{r-1}y^s\rp x \\
&+ \sum_{r\geq 0,\ s\geq0} b_{r,s}       x^ry^s - 
   \lp \sum_{r\geq 0,\ s\geq0} b_{r-1,s-1}   x^{r-1}y^{s-1}\rp xy
,
\end{align*}
where the third identity follows from (\ref{eqn:euler}) and the fourth
one follows from (\ref{eqn:a1}). If
$F= \sum_{r\geq 0,s\geq 0} b_{r,s}x^ry^s$ then
$E = xE + yE-xyE + F - xyF$ and hence
$(1-x-y-xy)E = (1-xy)F$. It follows from Subexample \ref{subexmp:misc}(1) that 
$F= \frac{1}{(1-x)^2(1-y)^2}$, therefore $E = \frac{1-xy}{(1-x)^3(1-y)^3}$.
\end{proof}
}}

The shortest path to compute the other Euler-Chow series is using the
full $\T = (\C^*)^2$-action on $\fl{0,1}{2}$, for it has finitely many
fixed points and orbits of dimension $1$, described in the following
lemma.
\begin{lemma}
\begin{enumerate}
\item The $\T$ action on $\fl{0,1}{2}$ has $6$ fixed points.
\item The $\T$ action has $9$ orbits of dimension $1$, divided
according to types as follows:
\begin{description}
\item[a] There are $3$ orbits of type $\scho{1;0,1}{2}$;
\item[b] There are $3$ orbits of type $\scho{0;0,2}{2}$;
\item[c] There are $3$ orbits of type $\scho{1;0,1}{2} + \scho{0;0,2}{2}$.
\end{description}
\end{enumerate}
\end{lemma}
\begin{proof}
Let $L_i = \{ z_i = 0 \}$, $i=0,1,2$, be the coordinate lines in $\Pp^2$ and
let $p_i, i=0,1,2$ denote the fixed points of the action in $\Pp^2$,
labeled such that $p_i\not\in L_i, i=0,1,2$. The orbits of the
first type consist of elements of the form $(l,L_i)$, $i=0,1,2$,
with $l\in L_i$.
The orbits of second type consist of elements of the form 
$(p_i,L), p_i \in L$. The orbits of third type consist
of elements of the form $(p,p_i\# p)$, where $p\in L_i$ and
$p_i\# p$ denotes the line determined by $p_i$ and $p$. 
\end{proof}

\begin{corollary}
\label{cor:f012all}
The $0$th and $1$st Euler-Chow series of $\fl{0,1}{2}$ are given
by the following generating functions:
$$\ecs{0}{\fl{0,1}{2}} = \frac{1}{(1-t)^6}$$
and
$$\ecs{1}{\fl{0,1}{2}} = \frac{1}{(1-r)^3(1-s)^3(1-rs)^3}\ ,$$
where $t$ is a variable associated to $\scho{0;0,1}{2}$ and 
$r,s$ are associated to $\scho{0;0,2}{2}$ and $\scho{1;1,2}{2}$,
respectively. 
\end{corollary}

One now can use the above information to compute the $p$-th Euler-Chow
series of the Grassmannian $\gr{1}{3}$, using the prescription of
Corollary \ref{cor:grass}.

In the next result we use the following association
 \{ variable \} $\leftrightarrow$ \{ Schubert class \}:
\begin{equation}
t \leftrightarrow \scho{0,1}{2}; \quad\quad
s \leftrightarrow \scho{0,2}{2}; \quad\quad
x \leftrightarrow \scho{0,3}{2}, \quad\quad
y \leftrightarrow \scho{1,2}{2}; \quad\quad
z \leftrightarrow \scho{1,3}{2}.
\end{equation}

\begin{proposition}
\label{prop:grass}
The Euler-Chow series of the Grassmannian $\gr{1}{3}$ are given by
the following generating functions:
\begin{align*}
\ecs{0}{\gr{1}{3}}  &= \frac{1}{(1-t)^6}, \quad 
& \ecs{1}{\gr{1}{3}}  &= \frac{1}{(1-s)^{12}}, \\
\ecs{2}{\gr{1}{3}}  &= \frac{1}{(1-x)^4(1-y)^4(1-xy)^3}, \quad
& \ecs{3}{\gr{1}{3}}  &= \frac{1+z}{(1-z)^5}, 
\end{align*}
where the latter coincides with
the Hilbert series, as explained in Example \ref{exmp:hfunction}.
\end{proposition}
\begin{proof}
It follows directly from Corollaries \ref{cor:grass} and
\ref{cor:tr-flag}, and the computations of the Euler-Chow functions of
the flag varieties $\fl{0,1}{2}$.
\end{proof}

\appendix
\section{Algebraic Constructions}
\label{app:A}

\subsection{Monoids with proper multiplication}
Throughout this discussion, 
all topological spaces considered are ``well behaved''
such as finite or countable CW-complexes.

Here we deal with {\bf abelian topological monoids}, although
we place special emphasis on discrete ones. However, the topological
approach yields a more unified perspective while retaining 
the origins of the subject. The category of all abelian 
topological monoids and  monoid morphisms will be denoted by \atm.

\begin{definition}
We say that $M \in \atm $ is a {\bf monoid with proper multiplication}
if the multiplication map is a {\bf proper} map; 
cf. Bourbaki \cite{bou-topI}.
We denote by \atmp\ the full subcategory of \atm\ consisting of
all abelian topological monoids with proper multiplication.
\end{definition}
\begin{example}
\label{1.2}
Any monoid $M \in \atm$ with finite multiplication, cf. Section
\ref{sec:prelim}, is an element of \atmp.
In particular, so are all free monoids.
\end{example}

We now discuss some functorial properties of the category \atmp.

\begin{proposition}
\label{pull-backs}\   
\begin{enumerate} 
\item If $\Psi : M \longrightarrow N$ is a {\bf proper} and {\bf surjective} 
monoid morphism between abelian topological monoids, 
and $M \in \atmp$, then so does $N$.
\item If $\Psi : M \longrightarrow N$ is a {\bf proper} monoid morphism
and $N \in \atmp$, then
so does $M$. In particular, the category \atmp\ is closed under pull-backs
over proper morphisms and under closed inclusions.
\end{enumerate}
\end{proposition}

\begin{proof}\ 
{\bf 1.}\quad Consider the following commutative diagram 
$$
\begin{CD} 
M \times M     @>{\ast_{M}}>>       M \\
@V{\Psi\times\Psi}VV       @VV{\Psi}V  \\
N \times N @>>{\ast_{N}}> N   
\end{CD}
$$

and let $C \subset N$ be compact. Since $\Psi$ is surjective, one has
$ \ast^{-1}_{N} (C) \, = 
        \, \Psi \times \Psi \, \lp   \ast^{-1}_{M}
                 \lp  \Psi^{-1}(C)\rp  \rp  
$
which is then compact, by the properness of $\Psi$ and $\ast_M$.

{\bf 2.}\quad We denote by $Q$ the pull-back in the following 
commutative diagram
\begin{center}
\mbox{
\xymatrix{
M \times M \ar[ddr]_-{\ast_{M}}\ar@{^{(}->}[r] & 
        Q \ar[r]^{\pi_{2}} \ar[dd]^{\pi_{1}} 
                        & M \times M \ar[d]^{\Psi\times\Psi} \\
& & N \times N \ar[d]^{\ast_{N}} \\
 & M \ar[r]_{\Psi} & N }
}
\end{center}
Since $\Psi$ and $\ast_{N}$ are both proper maps, it follows that
$\pi_{1}$ is also a proper map. The map
$j: M \times M \longrightarrow Q$ gives the homomorphism of
$M \times M $ onto the graph of $\ast_{M}$ contained in 
$Q \subset M \times M \times M$. In other words 
$j(m,\, m^{\prime}) \, = \, (m \ast_{M} \, m^{\prime}, \, m, m^{\prime})
\in Q$, and $j$ is a closed inclusion. The properness of
$\ast_{M}$ then follows from that of $\pi_{1}$ and the identity
$\pi_{1} \circ j = \ast_M$.

\end{proof}

\begin{remark}
It follows from the previous proposition that if a
monoid $M \in \atm$ has a {\bf proper} augmentation 
$\Psi : M \longrightarrow \Z_{+}$, then $M$ also belongs to \atmp.
This applies, in particular, to the Chow monoids introduced in Section
\ref{sec:chow}. 
\end{remark}
Now let $R$ be an arbitrary commutative ring and let $M$ be an
abelian topological monoid, which may be assumed discrete.
\begin{definition}
\label{aoR}\  
\begin{enumerate}
\item An {\bf $M$-graded algebra over $R$} is an $R$-algebra $A$ which can
be written as $A = \bigoplus_{m{\in}M} A_m$, where each $A_m$ is a
{\Zplus}-graded $R$-module and the multiplication 
$\cdot : A \times A \longrightarrow A$ induces a pairing of
$\Z_+$-graded $R$-modules 
$$
\cdot : A_{m,r} \times 
A_{{m^{\prime}},{r^{\prime}}} \longrightarrow
A_{m\,{\ast_{M}}\,{m^{\prime}},r+{r^{\prime}}}.
$$
Here $A_m = \bigoplus_{r\in{\Zplus}} A_{m,r}$ denotes the {\Zplus}-grading.
Let \gar{M}\ be the category of all
$M$-graded algebras over $R$ and grading 
preserving $R$-algebra homomorphism.

\item One says that an $R$-algebra $A \in \gar{M}$ is 
{\bf $M$-finite} or a {\bf finite $M$-graded algebra} if  
each   $A_m$ is an $R$-module of finite type. 
The full subcategory of \gar{M}\ consisting of those 
algebras which are {\bf $M$-finite} is denoted
by \garf{M} .
\end{enumerate}
\end{definition}

The following example plays a major role in this paper, since it
includes the Chow monoids as a particular case.
\begin{example}
\label{1.6}
Let $M$ be a monoid with {\bf proper multiplication} with the property
that all of its (path) components are compact. It follows from 
Proposition ~\ref{pull-backs} that the (discrete) monoid
$N = \pi_0 (M)$ is also in \atmp , and hence
its multiplication $\ast_N$ (induced by $\ast_M$) has
finite fibers, cf.~\ref{1.2}. Now, given a commutative
ring $R$, the singular homology
$$
A \, = \, H_{\ast} (M, R)
$$
of $M$ with coefficients in $R$ together with its
Pontriagrin ring structure becomes an $N$-graded algebra
over $R$, since $A = \bigoplus_{n\in{N}} H_{\ast} (M_n , R)$, where
$M_n$ denotes the connected component associated to
$n \in N = \pi_0 (M)$. Furthermore, since each $M_n$ is a compact
CW-complex, each $H_{\ast} (M_n , R) \equiv A_n$ is a finite
module over $R$. In other words, $A = H_{\ast} (M, R)$ is an
element of $\garf{N} = \gar{{\pi_0}(M)}$; cf. Definition ~\ref{aoR}. 
\end{example}

\begin{remark}
A typical example occurs when $M$ is an abelian 
topological monoid which comes
with a proper augmentation $\Psi : M \longrightarrow \Zplus$,
since $\Psi$ factors through $\pi_0 (M)$. See 
Proposition ~\ref{pull-backs}.
\end{remark}
We now discuss the \lq\lq change of grading\rq\rq\
behavior of our algebras under monoid 
morphisms $\Psi : M \longrightarrow N$.

\begin{definition}
\label{hola}
Let $A$ be an $M$-graded algebra over $R$ and $B$ be an $N$-graded algebra
over $R$. Given a monoid morphism $\Psi : M \longrightarrow N$, define
$\Psi^{\ast} B $ by
\begin{equation}
\Psi^{\ast} B \, = \, \bigoplus_{m\in{M}} \lp   \Psi^{\ast}B \rp   _m
\end{equation}
where $( \Psi^{\ast}B )_m = B_{\Psi(m)}$, and $\Psi_{\ast} A $ by
\begin{equation}
 \Psi_{\ast} A  \, = \,
        \bigoplus_{n\in{N}} \lp   \Psi_{\ast} A \rp  _n
\end{equation}
where $( \Psi_{\ast}A )_n = \bigoplus_{m\in{\Psi^{-1}(n)}} A_m$,
and $( \Psi_{\ast}A)_n = {0} \mbox{  if  } \Psi^{-1}(n) = \emptyset$.
Given homomorphisms $\varphi : A \longrightarrow A^{\prime}$ \,in \,
\gar{M}\ \, and \, $\eta : B \longrightarrow B^{\prime}$ \, in \,
\gar{N}, define:
\begin{equation}
\Psi^{\ast} (\eta) : 
\Psi^{\ast} B \longrightarrow \Psi^{\ast} B^{\prime} 
\mbox{ \, \,  by \, \,  } \Psi^{\ast} (\eta) (b) \, = \, \eta(b)
\end{equation}
and
\begin{equation}
\Psi_{\ast} (\varphi) : 
        \Psi_{\ast} A \longrightarrow \Psi_{\ast} A^{\prime}
\mbox{\, \, by \, \, } \Psi_{\ast} (\varphi)(a) \, = \, \varphi (a).
\end{equation}
\end{definition}
We have the following
 
\begin{proposition}
Let \, $\Psi : M \longrightarrow N$\, be a homomorphism of abelian
topological monoids.
\begin{enumerate}
\item The assignment of \, $ \Psi^{\ast} B \in \gar{M}$ \, to \, 
$B \in \gar{N}$, and of \, $\Psi^{\ast}\eta : \Psi^{\ast}B
\longrightarrow \Psi^{\ast}B^{\prime}$ \, to \, $\eta \in
\mbox{\bf Mor}_{\gar{N}} (B, B^{\prime})$ \, defines a
{\bf covariant} functor \, $\Psi^{\ast}: \gar{N} \longrightarrow
\gar{M}$. Furthermore, $\Psi^{\ast}$ preserves the algebras of 
finite type, in others words, $\Psi^{\ast}$ \, takes \,$\garf{N}$ \, 
to \garf{M}.

\item Similarly, \, $\Psi_{\ast}$  induces a {\bf covariant} functor
$\Psi_{\ast}: \gar{M} \longrightarrow \gar{N}$. If
$\Psi :M \longrightarrow N$ \, is a finite monoid morphism
(i.e. it has finite fibers), then $\Psi_{\ast}$ \, takes\,
\garf{M}\ \,into \,\garf{N}.
\end{enumerate}
\end{proposition}

\begin{proof}\ 
{\bf 1.}\quad We need to show the following: Given $B, B^{\prime} \in \gar{N}$
and a homomorphism $\varphi : B \longrightarrow B^{\prime}$ in
\gar{N}, one has:
\begin{description}
\item[a] $\Psi^{\ast} B \in \gar{M}$.
\item[b] $\Psi^{\ast}(\varphi) : 
        \Psi^{\ast}B \longrightarrow \Psi^{\ast}B^{\prime}$ is a morphism of
M-graded algebras over $R$ and $\Psi^{\ast}(1_B ) = 1_{\Psi^{\ast}B}$.
\item[c] Given $\varphi :B \longrightarrow B^{\prime}$\, and \,
$\varphi^{\prime} : B^{\prime} \longrightarrow B^{{\prime}{\prime}}$
one has \, $\Psi^{\ast}(\varphi^{\prime} \circ \varphi ) = 
\Psi^{\ast}(\varphi^{\prime}) \circ \Psi^{\ast}(\varphi)$.
\end{description}
By definition we have $(\Psi^{\ast}B)_m = B_{\Psi(m)}$ \, and the pairing
$B_{\Psi(m)} \times B_{\Psi(m^{\prime})} \longrightarrow
B_{{\Psi(m)}{\ast_{N}}{\Psi(m^{\prime})}} = 
                B_{\Psi(m{\ast_M}{m^{\prime}})}$ of $R$-modules, thus
one obtain a pairing
$
\lp   \Psi^{\ast} B \rp  _m \times  
\lp   \Psi^{\ast} B \rp  _{m^{\prime}} \longrightarrow
\lp   \Psi^{\ast} B \rp  _{m{\ast_M}{m^{\prime}}}
$
as required. Furthermore, if each $B_n$ is a {\bf finite}
$R$-module, for all $ n\in N$, then so is 
$( \Psi^{\ast} B )_{m}$ for all $m \in M$.
Now, given $\varphi : B \longrightarrow B^{\prime}$ an $R$-algebra
homomorphism such that $\varphi (B_n) \subset B^{\prime}_{n}$ then
if $b\in (\Psi^{\ast}B)_m = B_{\Psi(m)}$ and $ b^{\prime} \in
(\Psi^{\ast}B)_{m^{\prime}} =  B_{\Psi(m^{\prime})}$ one has that if
$bb^{\prime} \in B_{\Psi(m+m^{\prime})}$ then
$\Psi^{\ast}\varphi (bb^{\prime}) = \varphi (bb^{\prime}) = 
\varphi(b) \varphi(b^{\prime}) = 
\Psi^{\ast}\varphi (b) \cdot \Psi^{\ast}\varphi (b^{\prime})$,
and also $\Psi^{\ast}\varphi (b) = \varphi (b) \, \in \, 
B^{\prime}_{\Psi(m)} = (\Psi^{\ast}B^{\prime})_m$.
Therefore \, $\Psi^{\ast}\varphi$ \, is a homomorphism of $M$-graded
$R$-algebras. The remaining properties are proven in a similarly
trivial fashion.

{\bf 2.}\quad We need to show corresponding assertions for
$\Psi_{\ast} : \gar{M} \longrightarrow \gar{N}$. 
By definition, given $A \in \gar{M}$ one has 
$(\Psi_{\ast}A)_n = \bigoplus_{m\in{\Psi^{-1}(n)}} A_m$. Therefore,
given elements $\alpha \in (\Psi_{\ast}A)_n$ \, and \, 
$ \beta \in (\Psi_{\ast}A)_{n^{\prime}}$ \, one may assume
that there are $m \in \Psi^{-1}(n)$ \, and \, 
$m^{\prime} \in \Psi^{-1}(n^{\prime})$ such that 
$\alpha \in A_m$ \, and \, $ \beta \in A_{m^{\prime}}$.
In this case $\alpha \cdot \beta \in A_m \cdot A_{m^{\prime}}
\subset A_{m{\ast_M}m^{\prime}}$, \, and hence 
$\alpha \beta \in (\Psi_{\ast}A)_{n{\ast_N}n^{\prime}}$,
\, since $\Psi(m{\ast_M}m^{\prime}) = \Psi(m)\ast_N \Psi(m^{\prime})
= n{\ast_N}n^{\prime}$. \, In other words, \, $\Psi_{\ast}A$ \, is
an $N$-graded $R$-algebra. Now, given \, 
$\varphi : A \longrightarrow A^{\prime}$ \, and \,
$\varphi^{\prime} : A^{\prime} \longrightarrow A^{\prime\prime}$ \,
morphisms in \gar{M}\ one has 
$\Psi_{\ast}(\varphi^{\prime} \circ \varphi) =  
\Psi_{\ast}(\varphi^{\prime}) \circ \Psi_{\ast}(\varphi)$ \, by
definition, and if \,$\alpha \in (\Psi_{\ast}A)_n$ \, then
\, $\alpha \in A_m$ \, for some \,$m\in \Psi^{-1}(n)$\, (otherwise
$\alpha = 0$, by definition). Therefore \,$\varphi(\alpha) \in 
A^{\prime}_m \subset (\Psi_{\ast}A^{\prime})_n$, \, and
hence \,$ \Psi_{\ast}(\varphi)$ \, is a morphism in \, \gar{N}.
Finally, if $\Psi$ has finite fibers, then\,
$(\Psi_{\ast}A)_n = \bigoplus_{m\in\Psi^{-1}(n)} A_m$\, is a sum
over finitely many indices, showing that $\Psi$ then sends
\garf{M}\ into \garf{N}.
\end{proof}

\subsection{Invariants for $M$-graded algebras}

We now restrict our attention to elements in \atmp\ which have
{\bf finite multiplication}, and refer the reader to
Section \ref{sec:prelim} for the notation used here.
Our goal is to introduce the following invariants of $M$-graded
algebras over $R$.

\begin{definition}
\label{defn:hilbert}
Let $M \in \atmp$ be a monoid with finite multiplication and 
let $A \in \garf{M}$ be an finite $M$-graded algebra over $R$.
\begin{enumerate}
\item If $R$ is a principal ideal domain (PID) define the {\bf Hilbert 
$M$-function} (or Hilbert $M$-series)
$P_A (t) \in \Z [t]^M$ of $A$ to be the function $P_A (t):
M \longrightarrow \Z [t]$ which sends $m\in M$ to
$$
P_A (t)(m) \, = \, \sum_{k\in{\Zplus}} \lp  {\rank}_R \, A_{m,k}\rp   t^k .
$$
\item Under the same hypothesis, define the {\bf Euler $M$-function}
(or Euler $M$-series)
$E_A \in \Z^{M}$ of $A$ to be the image of $P_A (t)$ in $\Z^M$ under 
the \lq\lq evaluation homomorphism\rq\rq\ at -1:
$ (e_{-1})_{\ast}: \Z [t]^M \longrightarrow \Z^M$ 
cf. Proposition \ref{13} (3). In other words, given $m \in M$,
$$
E_A (m) \, = \, P_A (-1)(m) \, = \, 
        \sum_{k\in{\Zplus}} (-1)^k \rank \lp  A_{m,k}\rp  .
$$
\end{enumerate}
\end{definition}
We have the following
\begin{proposition}
Let $\Psi : M \longrightarrow N$ be a finite morphism of monoids with
finite multiplication, and let $A \in \garf{M}$ be a finite $M$-graded algebra
over $R$. Then
$$
P_{\Psi_\sharp{A}}(t) \, = \, \Psi_{\sharp} \lp  P_A (t)\rp  .
$$
Similarly, if $B \in \garf{N}$ then
$$
P_{\Psi^\sharp{B}}(t) \, = \, \Psi^{\sharp} \lp  P_B (t)\rp  .
$$
\end{proposition}

\begin{proof}
Given $n \in N$, by definition one has
\begin{align*}
P_{\Psi_{\sharp}A} (t)(n) 
& = \sum_{k\in{\Zplus}} \rank_R \lp  \lp  \Psi_{\sharp}A\rp  _{n,k}\rp   t^k \;
 = \; \sum_{k\in{\Zplus}}\rank_R 
        \lp  \sum_{m\in\Psi^{-1}(n)} A_{m,k}\rp   t^k\\
& = \sum_{m\in\Psi^{-1}(n)} \, \, 
        \sum_{k\in{\Zplus}} \, \lp   \rank_R A_{m,k}\rp   t^k 
\; = \; \sum_{m\in\Psi^{-1}(n)} \, P_A (t)\, (m)\\
&= \lp  \Psi_{\sharp} \lp  P_A (t)\rp  \rp  (n).
\end{align*}
Similarly,
\begin{align*}
\lp  \Psi^{\sharp} \lp  P_B (t)\rp  \rp  (m)
& = P_B (t) \lp  \Psi(m)\rp  
\; =\; \sum_{k\in\Zplus} \rank_R \lp  B_{\Psi{(m)},k}\rp  \, t^k\\
& = \sum_{k\in\Zplus} \rank_R \lp   (\Psi^{\sharp}B)_{m,k}\rp   \, t^k
\; = \; P_{\Psi^{\sharp}B} (t).
\end{align*}
\end{proof}

\begin{remark}
\label{rem:ecf}
Note that if $M$ is the Chow monoid $\cvp{p}{X}$ of a projective
variety $X$, then its Pontrjagin ring $H_*(M;\Z)$ is a finite
$\ax{p}{X}$-graded algebra whose Euler $\ax{p}{X}$-function is precisely
the $p$-th Euler-Chow function $\ecs{p}{X}$ of $X$.
\end{remark}



\providecommand{\bysame}{\leavevmode\hbox to3em{\hrulefill}\thinspace}


\end{document}